\documentclass[12pt,a4paper]{article}

\usepackage{xcolor}

\usepackage{amsfonts}
\usepackage{amsmath}
\usepackage{amsbsy}
\usepackage{theorem}
\usepackage{graphicx}

\newtheorem{theorem}{Theorem}
\newtheorem{proposition}{Proposition}
\newtheorem{lemma}{Lemma}
\newtheorem{corollary}{Corollary}
\newtheorem{remark}{Remark}

\begin{document}
\sf
\title{\bfseries Long range dependence spectral  testing from  discretely observed  functional time series in manifolds
}
\date{}

 \maketitle
 \author{\hspace*{-1cm} M.D. Ruiz-Medina$^{1}$ and Rosa M. Crujeiras$^{2}$\\
$^{1}$ University of Granada\\
$^{2}$ University of Santiago de Compostela}

 \begin{abstract}
A specification test  for  long range dependence (LRD) in  functional time series in manifolds has been formulated in \cite{Ruiz-MedinaCrujeiras24} in the spectral domain for fully observed functional data. The asymptotic Gaussian distribution of the proposed test statistics, based on the weighted periodogram operator, under the null hypothesis, and the consistency of the test have been derived. In this paper, we analyze the asymptotic properties of this  LRD  spectral testing procedure, when the  functional data are  contaminated, and discretely observed through random uniform spatial sampling.
\end{abstract}
\medskip


\medskip

\newpage

\subsection*{Graphical abstract}
\begin{figure}[!h]
\begin{center}
\includegraphics[width=16cm,height=5cm]{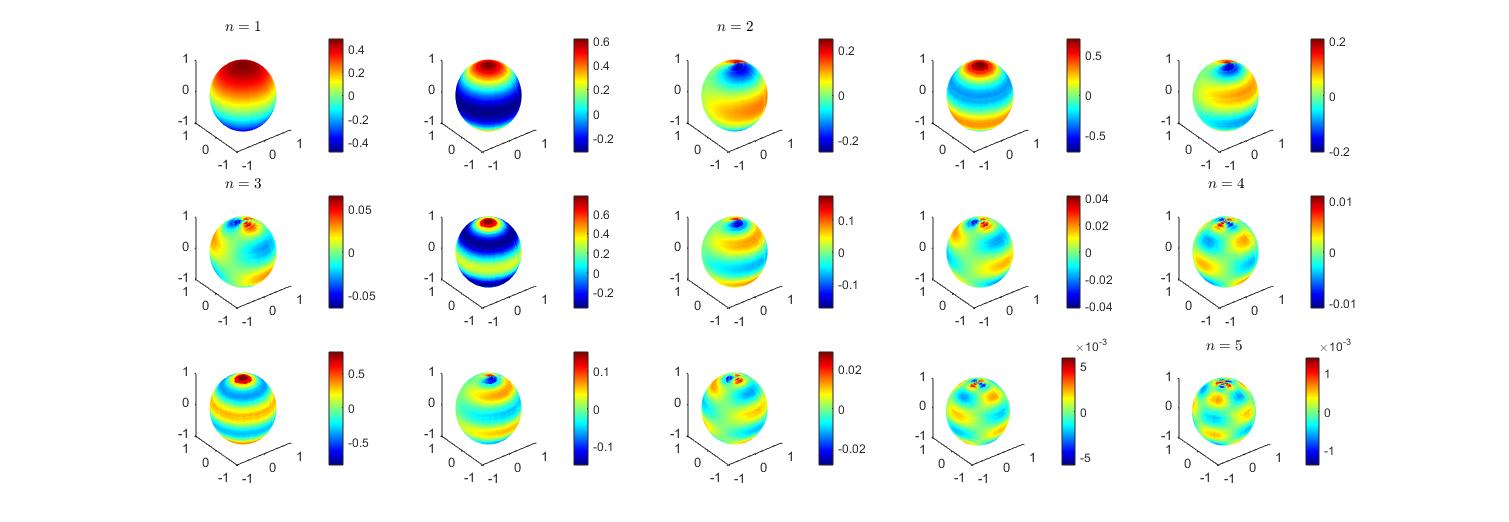}
\includegraphics[width=16cm,height=3.5cm]{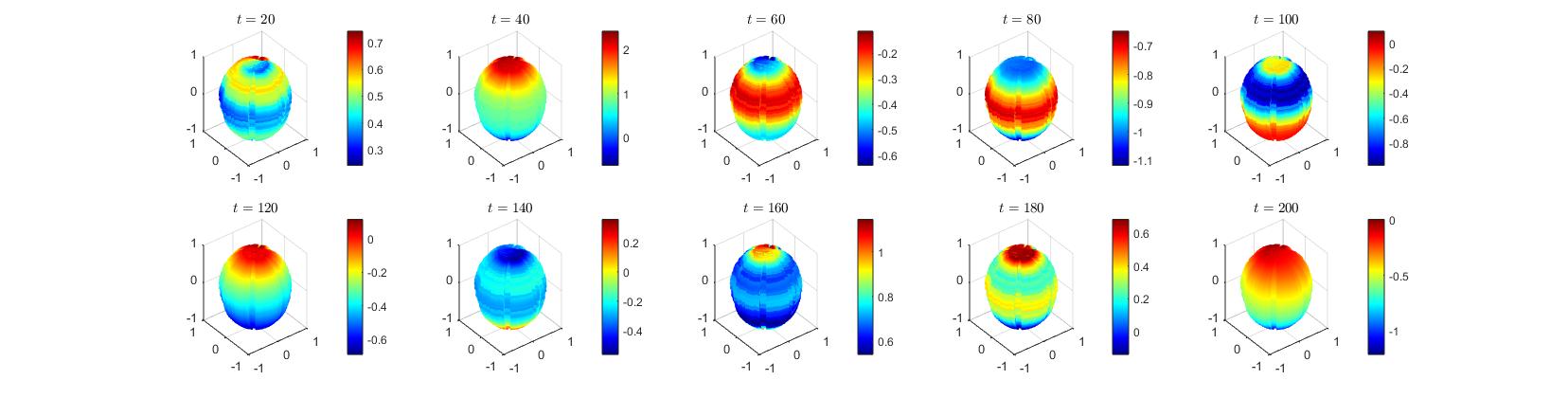}
\includegraphics[width=16cm,height=3.5cm]{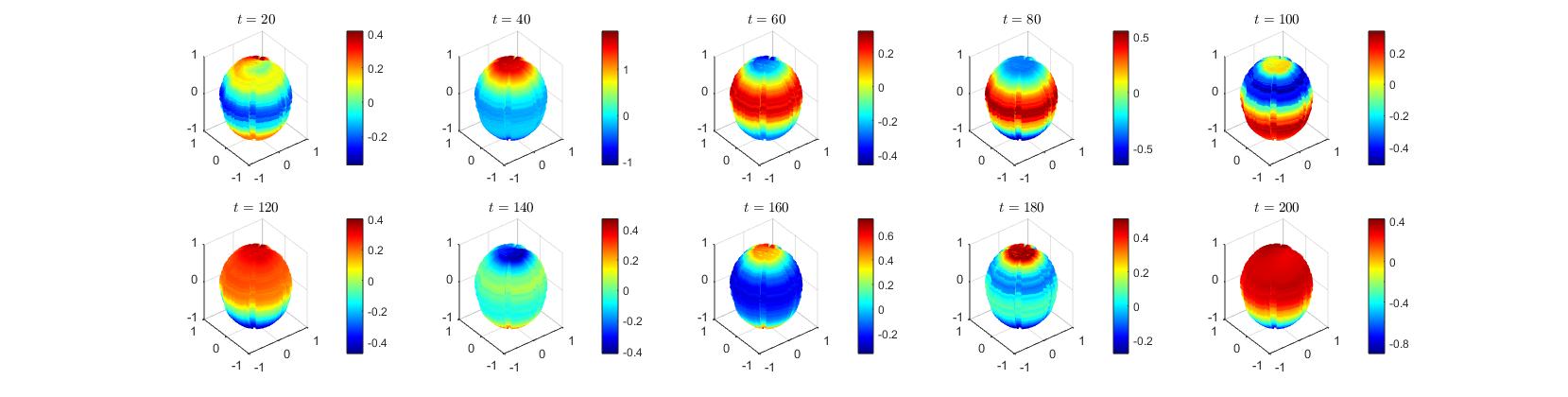}
\end{center}
\caption{The truncated random spherical harmonics sieve basis of dimension $k(T)=15$ (three first lines at the top), the contaminated discretely observed data (next two lines at the center), and its nonparametric series  Least-Squares (LS)  reconstruction (last two lines at the bottom)}
\end{figure}

\clearpage

\subsection*{Highlights}
\begin{itemize}
\item Research highlight 1. Testing Long Range Dependence (LRD) in functional time series constitutes a challenging topic where several problems remain open. This testing problem has been recently  studied in Ruiz-Medina and Crujeiras (2025).  Its asymptotic analysis has been achieved in the double spectral domain for fully observed    functional time series in manifolds. This paper  addresses the case of irregularly discretely observed and contaminated functional data correlated in time.
\item Research highlight 2. Low values of the sieve parameter are considered for spatial sparsity,  under smooth  functional data,  when the strongest persistent in time is displayed at coarser spatial scales. A balance between dimensionality in space and time, and spatial smoothness and memory is required to ensure the asymptotic results hold.
\end{itemize}

\newpage

\section{Introduction}
\label{sec:1}
Hypothesis testing in the spectral domain for stationary functional time series has  mainly been developed in the weak--dependent case.
 In
\cite{Panaretos13}, the asymptotic properties of the weighted periodogran operator are obtained
for stationary functional time series, under  Short Range Dependence (SRD),  assuming  suitable summability conditions of the $L^{2}$ norms of the cumulant operators. Specifically,  the asymptotic   Gaussian distribution of the functional discrete Fourier transform (fDFT) and weighted periodogram operator are derived.
Bias asymptotic analysis  in $L^{2}$  norm of the empirical spectral cumulant  operators is addressed. In particular, the consistency of the weighted periodogram operator  in the integrated mean square error sense is established.     A spectral version of covariance tests,  based on the Hilbert--Schmidt operator norm,  to compare the second--order structure of two functional time series in terms of their spectral density operator families, is formulated in \cite{TP16} (see also \cite{Boenteetal18};
\cite{PS16}; \cite{Pigolietal14}).

Beyond the weak--dependence and stationary assumptions, LRD functional time series analysis is developed in \cite{LIROB19}  via the pure point  spectral decomposition of the long-run covariance operator. This paper opens new research lines regarding the consistent estimation of the dimension and the orthonormal
functions  spanning the dominant subspace, where the projected curve
process displays the largest dependence range. The results are applied to the statistical analysis of fractionally integrated functional autoregressive moving averages processes.

 In the stationary case, a characterization of LRD in the spectral domain for functional time series is derived  in \cite{RuizMedina2022}.
The asymptotic unbiasedness of the integrated periodogram operator in the Hilbert--Schmidt operator norm is obtained.
The spectral density operator family is assumed to be diagonalized in terms of a common resolution of the identity (spectral kernel). The LRD operator also admits a weak--sense diagonalization in terms of this spectral kernel, and the corresponding pure point spectrum  is parameterized.
Under a Gaussian scenario, the weak--consistent minimum contrast  estimation of the LRD operator is derived  in a   semiparametric  spectral  framework.  The results derived in \cite{RuizMedina2022} can be applied to stationary functional time series in  connected and compact two point homogeneous spaces, under  the  invariance assumption with respect to the group of isometries of the manifold (see \cite{OvalleRM24}).

A statistical hypothesis test, based on the weighted periodogram operator,  for detecting LRD in functional time series in manifolds,
  is formulated  in the double spectral  domain in \cite{Ruiz-MedinaCrujeiras24}. The  asymptotic normality of the proposed test statistics under the null hypothesis is derived,  from an extended version of the results  in \cite{Panaretos13} for SRD stationary functional time series.
   Bias asymptotic analysis  of the empirical spectral  cumulant operators   in \cite{Panaretos13} under SRD is also extended in  \cite{Ruiz-MedinaCrujeiras24} to the LRD functional time series context. The consistency of the integrated weighted periodogram operator under LRD  is then obtained, leading  to   the almost surely divergence of the Hilbert--Schmidt operator norm of the proposed test statistics under the alternative, and  yielding
the consistency of the test.

The results derived in all the above cited references on spectral analysis of functional time series and, in particular, on  LRD analysis  are obtained under the umbrella of fully observed functional data. This paper addresses    the case of (irregularly) discretely  observed functional time series  contaminated by additive noise. To implement our LRD testing procedure, nonparametric series least squares regression, based on   random uniform manifold sampling, is first applied in the   reconstruction of the functional data values. The Gaussian asymptotic distribution of the resulting plug--in test statistics under the null hypothesis is derived.  $L^{2}$--convergence rates are obtained in the  bias asymptotic analysis of the integrated second and fourth order plug-in empirical spectral  cumulant  operators under LRD. Consistency of the test is also derived  from  the  consistency of the integrated  plug--in weighted periodogram operator under LRD.  This asymptotic analysis is achieved from the suitable  interaction between the   increasing rate of the  functional sample size, the spatial sampling frequency, and  the sieve parameter, as well as the decreasing rate  of the bandwidth parameter. The  local  smoothness of the functional data,  and the memory in time of the underlying model also play a crucial role in this analysis, as commented in the Concluding Remarks section. The simulation study undertaken in Section  \ref{simulationstudy}  illustrates consistency, and the   finite sample performance of the  LRD testing procedure, by computing empirical test sizes and powers under sparse spatial data scenarios.

The summary of the contents of the paper is the following. The methodological approach, including the formulation of the testing procedure and the observation model considered, are introduced in Section \ref{sec2}. The mean-square consistency  of the nonparametric series least squares estimator of the functional values of the data is derived in Section \ref{dodata}. The asymptotic properties of the plug--in test statistics  under the null and the alternative hypotheses  are obtained in Section \ref{asymptotics}. The practical illustration of the methodological approach proposed is provided in the simulation study in Section \ref{simulationstudy}.
Section \ref{cr} establishes some concluding remarks.  Finally, the functional spectral background, required in the formulation of the LRD testing procedure, is introduced in  \ref{SBLRDF}. The proofs of the results are given in   \ref{supproofres}.

\section{Methodological approach}
\label{sec2}
 In what follows, we will introduce the notation to be used along the paper.
 $L^{2}\left(\mathbb{M}_{d},d\nu , \mathbb{R}\right)$ and $L^{2}(\mathbb{M}_{d},d\nu; \mathbb{C})$ respectively denote the space   of real-- and complex-- valued  square integrable functions  on a Riemannian manifold $\mathbb{M}_{d},$   given by a connected and compact  two--point homogeneous space embedded into $\mathbb{R}^{d+1}.$  The topological dimension of $\mathbb{M}_{d}$ is
  $d,$ and $d\nu$  denotes the normalized Riemannian measure on $\mathbb{M}_{d}$ (see, e.g.,  \cite{MaMalyarenko}).
 Let $X=\left\{X_{t},\ t\in \mathbb{Z}\right\}$  be an $L^{2}\left(\mathbb{M}_{d},d\nu , \mathbb{R}\right)$--valued correlated sequence,
  defined on the basic probability space  $(\Omega ,\mathcal{Q}, \mathcal{P}),$
 which is assumed to be
   strictly stationary (in time and space)  with  zero--mean.
   \subsection{Formulation of the testing procedure}
We adopt the approach introduced in \cite{Ruiz-MedinaCrujeiras24},  in the formulation of the LRD testing procedure in the functional spectral domain (see  \ref{SBLRDF} for more details).

 We are interested on testing SRD against LRD.  It is well-known that  SRD in real--valued stationary  time series is characterized by the absolutely
summability of the covariance function. This notion is extended to the infinite--dimensional framework in terms of the summability of the sequence of nuclear norms of the elements of the  covariance operator family. That is, SRD is understood in the following sense:
\begin{equation}\sum_{\tau \in \mathbb{Z}}\|\mathcal{R}_{\tau}\|_{L^{1}(L^{2}(\mathbb{M}_{d},d\nu , \mathbb{R}))}=\sum_{\tau \in \mathbb{Z}}\sum_{n\in \mathbb{N}_{0}}\Gamma (n,d)\left|\int_{-\pi}^{\pi}\exp\left(i\omega \tau \right)f_{n}(\omega )d\omega \right|<\infty,
\label{srd}
\end{equation}
\noindent where the  equality follows from (\ref{sdo2}) under the invariance property of covariance and spectral density kernels.  Here, $L^{1}(L^{2}(\mathbb{M}_{d},d\nu , \mathbb{R}))$ denotes the space of trace or nuclear operators on $L^{2}(\mathbb{M}_{d},d\nu , \mathbb{R}).$ Series $X$ is said to display LRD
when this summability condition fails.
\begin{remark}
The summability condition (\ref{srd})  can  be rewritten   as \linebreak  $\sum_{\tau \in \mathbb{Z}}\int_{\mathbb{M}_{d}}r_{\tau}(x,x)d\nu(x),$ being approximated, in the case of discretely observed functional data,  via Monte Carlo numerical integration as  \linebreak $\sum_{\tau \in \mathbb{Z}}\frac{1}{M}\sum_{m=1}^{M}r_{\tau}(x_{m},x_{m}),$
\noindent where $x_{1},\dots, x_{M},$ denote the observed or generated  values of  a sample  of size  $M$ of independent and identically uniform distributed  random variables in $\mathbb{M}_{d}.$
Equivalently, the  underlying   time-varying local spatial  regularity order, reflected on the rate of convergence to zero of the spatial pure point spectra of the elements of the  covariance operator family,  is summable in time, leading to the  summability of  $\sum_{\tau \in \mathbb{Z}}\frac{1}{M}\sum_{m=1}^{M}r_{\tau}(x_{m},x_{m}).$ Hence,  in the compact spatial domain case, SRD   also provides information about local spatial smoothness.

The situation is different in the case of noncompact spatial domain, when the SRD condition
$$\sum_{\tau \in \mathbb{Z}}\int_{\mathbb{R}^{d+1}}r_{\tau}(\|x\|)dx,$$\noindent  is assumed under  homogeneity and isotropic in space.

\end{remark}

 From equation (\ref{LRD}), the positive  function sequence $\{f_{n}(\omega ),\ \omega \in [-\pi,\pi],\ n\in \mathbb{N}_{0}\}$ in  (\ref{sdo2}) satisfies:
  \begin{eqnarray}  f_{n}(\omega )&=&\frac{M_{n}(\omega )}{|\omega |^{\alpha (n)}},\quad \omega \in [-\pi,\pi],\ n\in \mathbb{N}_{0}.
  \label{eqscm1}
 \end{eqnarray}

  The following testing problem is then  considered:
 \begin{eqnarray}
 && H_{0}: \ f_{n}(\omega )=M_{n}(\omega ), \  \omega \in [-\pi,\pi],\ \forall n\in \mathbb{N}_{0}\label{altfhp0}\\
 && H_{1}: \ f_{n}(\omega )=M_{n}(\omega )\left|\omega\right|^{-\alpha (n)},\ \omega \in [-\pi,\pi],\ \forall   n\in \mathbb{N}_{0}.
 \label{reatlrd}
 \end{eqnarray}

\noindent    We adopt the design of the test statistics $\mathcal{S}_{B_{T}}$  proposed in  \cite{Ruiz-MedinaCrujeiras24},  capturing the singularities at zero frequency for different manifold resolution levels under $H_{1}.$ Specifically, $\mathcal{S}_{B_{T}}$ is given by
\begin{equation}\mathcal{S}_{B_{T}}=\sqrt{B_{T}T}\int_{[-\sqrt{B_{T}}/2, \sqrt{B_{T}}/2]}\widehat{\mathcal{F}}_{\omega }^{(T)}
\frac{d\omega}{\sqrt{B_{T}}}= \int_{-\pi}^{\pi}\frac{\sqrt{B_{T}T}\widehat{\mathcal{F}}_{\omega }^{(T)}\mathbb{I}_{[-\sqrt{B_{T}}/2, \sqrt{B_{T}}/2]}(\omega )}{\sqrt{B_{T}}}d\omega,\label{oslrdtintro}
\end{equation}

\noindent where the kernel of the weighted periodogram operator  $\widehat{\mathcal{F}}_{\omega }^{(T)}$ has been introduced in equation (\ref{enp}), with  $B_{T}$ being   the bandwidth parameter. The indicator function on the interval $[-\sqrt{B_{T}}/2, \sqrt{B_{T}}/2]$ is denoted by  $\mathbb{I}_{[-\sqrt{B_{T}}/2, \sqrt{B_{T}}/2]}(\omega ),$ for $\omega \in [-\pi,\pi].$  Note that,  as $T\to \infty,$  $$\int_{-\pi}^{\pi}\frac{\mathbb{I}_{[-\sqrt{B_{T}}/2, \sqrt{B_{T}}/2]}(\omega )}{\sqrt{B_{T}}}h(\omega  )d\omega \to \int_{-\pi}^{\pi}\delta (0-\omega )h(\omega  )d\omega = h(0),$$ \noindent for every $h\in L^{2}([-\pi ,\pi])$  (see \cite{Gelfand64}).   Here,  $\delta (0-\omega )$ denotes  the  Dirac Delta distribution
at zero frequency. Hence, we will adopt the notation $\delta_{T}(0-\omega )=\frac{\mathbb{I}_{[-\sqrt{B_{T}}/2, \sqrt{B_{T}}/2]}(\omega )}{\sqrt{B_{T}}}$ representing a truncated Dirac Delta distribution.

\begin{remark}
The strong consistency of $\mathcal{S}_{B_{T}}$ is obtained in \cite{Ruiz-MedinaCrujeiras24}  under fully observed functional data, exploiting the asymptotic behavior of $\delta_{T}(0-\omega ),$ and of the weighting function  $W^{(T)}$  in (\ref{enp}), as $T\to \infty.$  Additional conditions are derived  in Theorem \ref{th2} below to ensure strong consistency of $\mathcal{S}_{B_{T}}$ when functional data are discretely observed.
  \end{remark}

\subsection{Nonparametric series  Least-Square (LS) reconstruction of functional data}
\label{nslsfd}
In this section, we apply nonparametric series least-squares, in the  reconstruction over $\mathbb{M}_{d}$ of our discretely observed functional data set. The mean--square consistency    is also derived.

Let us consider the time--varying  Karhunen--Lo\'eve expansion of $X=\left\{X_{t},\ t\in \mathbb{Z}\right\}$ in $L^{2}\left(\Omega, \mathbb{M}_{d},[0,T],\mathcal{P}\otimes d\nu\otimes dt\right)$ (see, e.g., \cite{MarinucciRV}):

\begin{eqnarray}
X_{t}(z,\xi) &=&\sum_{n\in \mathbb{N}_{0}}\sum_{j=1}^{\Gamma (n,d)}\left\langle X_{t}(\cdot, \xi),S_{n,j}(\cdot) \right\rangle_{L^{2}\left(\mathbb{M}_{d},d\nu , \mathbb{R}\right)}S_{n,j}(z)\nonumber\\
&=& \sum_{n\in \mathbb{N}_{0}}\sum_{j=1}^{\Gamma (n,d)}X_{nj}(t,\xi)S_{n,j}(z),\quad \xi\in \Omega ,\ z\in \mathbb{M}_{d},\ t\in [0,T].\label{eqkshexpansion}
\end{eqnarray}
\noindent  Here, from equation (\ref{sdo2}),  for $n\in \mathbb{N}_{0},$  \begin{eqnarray}&& \int_{\Omega }X_{nj}(t,\xi)X_{nj}(t+s,\xi)\mathcal{P}(d\xi) =E[X_{nj}(t,\xi)X_{nj}(t+s,\xi)] \nonumber\\
&&= E\left[\left\langle X_{t}(\cdot, \xi),S_{n,j}(\cdot) \right\rangle_{L^{2}\left(\mathbb{M}_{d},d\nu , \mathbb{R}\right)}\left\langle X_{t+s}(\cdot, \xi),S_{n,j}(\cdot) \right\rangle_{L^{2}\left(\mathbb{M}_{d},d\nu , \mathbb{R}\right)}\right]\nonumber\\
&&\int_{\mathbb{M}_{d}\times \mathbb{M}_{d}}r_{s}\left(x,y\right)S_{n,j}(x)S_{n,j}(y)d\nu(x)d\nu(y)
\nonumber\\
 && =\int_{-\pi}^{\pi}\exp(is\omega )f_{n}(\omega)d\omega ,\quad j=1,\dots,\Gamma(n,d).\end{eqnarray}

\noindent For each $t\in [0,T],$
\begin{eqnarray}
&&\int_{\Omega \times \mathbb{M}_{d}}\left|X_{t}(z,\xi) -\sum_{n=1}^{N(T)}\sum_{j=1}^{\Gamma (n,d)}X_{nj}(t,\xi)S_{n,j}(z)\right|^{2}d\nu(z)\mathcal{P}(d\xi)\nonumber\\
&& =E\left[\int_{\mathbb{M}_{d}}\left|X_{t}(z,\xi) -\sum_{n=1}^{N(T)}\sum_{j=1}^{\Gamma (n,d)}X_{nj}(t,\xi)S_{n,j}(z)\right|^{2}d\nu(z)\right]
 = \mathcal{O}\left([k(T)]^{-\frac{2s}{d}}\right),\nonumber\\
\label{eqkshexpansionbb}
\end{eqnarray}
\noindent where we have applied that, for a function  in the  fractional Sobolev $H^{s}(\mathbb{M}_{d})$ of order $s>d/2,$   the convergence rate, in $L^{2}(\mathbb{M}_{d}, d\nu)$ norm, of its finite-dimensional approximation, by projection onto the subspace generated by $\sum_{l=0}^{N(T)}\Gamma (l,d)=k(T)$ elements of  the orthonormal basis  $\{S_{n,j}^{d}, \ j=1,\dots,$ \linebreak  $\Gamma (n,d),\ n\in \mathbb{N}_{0}\}$ of $L^{2}(\mathbb{M}_{d}, d\nu),$  is given by the asymptotic order
$\mathcal{O}\left([k(T)]^{-\frac{2s}{d}}\right)$ (when the square of the norm is considered).
That  is the case of functions in the Reproducing Kernel Hilbert Space (RKHS), generated  by the kernel of the autocovariance operator $\mathcal{R}_{0}$ of $X=\left\{X_{t},\ t\in \mathbb{Z}\right\}$ in our framework (see also condition (ii) below). Here, since $k(T)$ denotes  the sieve parameter depending on the functional sample size $T,$  we also adopt the notation $N(k(T))$ to indicate dependence of $N(T)$  on   $k(T).$
 Note that, from Chebychev inequality,  applying Borel Cantelli Lemma, a.s. convergence in (\ref {eqkshexpansionbb}) also holds for $s>d/2,$ from Embedding Theorems of fractional Sobolev spaces into H\"older  spaces (see, e.g.,  \cite{Triebel78, Triebel97}).

\section{Theoretical results for discretely observed functional data}
\label{dodata}
The formulation of the observation  model, constituted by time-correlated, contaminated, discretely and irregularly observed functional data, and the  assumptions required  are first introduced.  Theorem \ref{th1} is then obtained, providing the mean-square rate of convergence of the consistent nonparametric series  LS functional spatial  interpolator.

\subsection{Discretely observed functional data}

   Let us  consider the following observation model:
   \begin{equation}Y_{t}(Z_{i})=X_{t}(Z_{i})+\varepsilon_{i,t}, \ i=1,\dots, M(T),\ t=0,\dots,T-1,\label{reg}
    \end{equation}\noindent where  $\varepsilon_{i,t},$ $i=1,\dots,M(T),$ $t=0,\dots,T-1,$ are independent and identically distributed random variables, and  $E\left[\varepsilon_{i,t}/Z_{i}\right]=0=E\left[\varepsilon_{i,t}\right],$
   $ i=1,\dots, M(T),$ $t=0,\dots,T-1.$ Here, $M(T)$ is the number of spatial sampling locations in $\mathbb{M}_{d}$ when the functional sample size is $T.$ Assume that, for every $T\geq 2,$
  $Z_{1},\dots, Z_{M(T)}$ are independent and uniform distributed in $\mathbb{M}_{d}.$   From the  sample
   $\left(Y_{t}(Z_{i}), Z_{i}\right),$  $i=1,\dots,M(T),$ $t=0,\dots, T-1,$ obeying  (\ref{reg}),   we compute the nonparametric series LS estimator
  \begin{equation}\widehat{\mathbf{X}}_{t,k(T),M(T)}(z)=b_{k(T)}^{\prime }(z)(B_{k(T),M(T)}^{\prime }B_{k(T),M(T)})^{-}B_{k(T),M(T)}^{\prime }\mathbf{Y}_{t},\label{LSFC}
  \end{equation}
 \noindent where  \begin{equation}b_{k(T)}^{\prime }(z)=\left(S_{n,j}^{d}(z), \ j=1,\dots, \Gamma (n,d),\ n=0,\dots, N(k(T))\right),\label{eqvectsb}
 \end{equation} \noindent is a row $1\times k(T)$ vector, for each $z\in \mathbb{M}_{d},$
 and   $B_{k(T),M(T)}=$\linebreak $\left(b_{k(T)}(Z_{1}),\dots,b_{k(T)}(Z_{M(T)})\right)^{\prime }$ is a $M(T)\times k(T)$ matrix.
  The exponent $^{-}$ denotes the Moore--Penrose generalized inverse, and
   $\mathbf{Y}_{t}=\left(Y_{t}(Z_{1}),\dots,Y_{t}(Z_{M(T)})\right)^{\prime },$  $t=0,\dots, T-1.$
  \subsection{Assumptions}

In the derivation of the main result of this section, Theorem \ref{th1}, the following assumptions are made:
  \begin{itemize}
\item[(i)] $\sup_{z\in \mathbb{M}_{d}}E\left[\left|\varepsilon_{i,t}\right|^{2}/Z_{i}=z\right]<\infty,$ $i=1,\dots, M(T),$ $t=0,\dots,T-1.$
          \item[(ii)] For $t\in \mathbb{Z},$  $\mathcal{P}\left( X_{t}\in H^{s}(\mathbb{M}_{d})\right)=1,$  with $H^{s}(\mathbb{M}_{d})$ denoting the fractional Sobolev  of order $s>d/2$ of functions with compact support contained in $\mathbb{M}_{d}.$
          \item[(iii)] $[k(T)]^{2}/M(T)\to 0,$ $T\to \infty.$
         \end{itemize}
\begin{remark}\label{condiii}  Although Theorem \ref{th1} below is obtained under  (iii),  the  next section analyzes  the special case of condition (iii) given by     $k(T)=T^{\widetilde{\alpha }}$ and  $M(T)=T^{\gamma },$ with $\gamma -2\widetilde{\alpha }>0$  (see  Remark \ref{remconiiiapp} below). The simulation study undertaken  in Section \ref{simulationstudy}  goes beyond this particular scenario of (iii).
\end{remark}
     Mean-square consistency of Monte Carlo numerical integration, based on random  uniform sampling in $\mathbb{M}_{d},$ is applied to obtain the following rate of convergence \begin{eqnarray}
     &&\hspace*{-0.8cm} E\left[\left\|\frac{B_{k(T),M(T)}^{\prime }B_{k(T),M(T)}}{M(T)}-I_{k(T)\times k(T)}\right\|^{2}\right]=
     \mathcal{O}\left( \frac{[k(T)]^{2}}{M(T)}\right),\ T\to \infty,\label{asidentity5}
     \end{eqnarray}
      \noindent of the discrete approximation $\left\|\frac{B_{k(T),M(T)}^{\prime }B_{k(T),M(T)}}{M(T)}-I_{k(T)\times k(T)}\right\|^{2}$ to the \linebreak $L^{2}\left(\mathbb{M}_{d},d\nu , \mathbb{R}\right)$ inner products between the basis elements \linebreak $S_{n,j}^{d}(z), \ j=1,\dots, \Gamma (n,d),\ n=0,\dots, N(k(T)),$ for each  $T\geq 2.$ Here,
        $ \|\cdot \|$ denotes  the Frobenius (or Euclidean) norm. Thus, the equivalence of the empirical and the theoretical  $L^{2}$ norms holds over the linear subspace generated by $$\left(S_{n,j}^{d}(z),  \ j=1,\dots, \Gamma (n,d),\ n=0,\dots, N(k(T))\right).$$

\begin{theorem}
\label{th1}
Under conditions (i)--(iii),  for $t=0,\dots,T-1,$
  \begin{eqnarray}&&\widehat{\mathbf{X}}_{t,k(T),M(T)}(x)=X_{t}(x)+\mathcal{O}\left(\sqrt{\frac{k(T)}{M(T)}}\right)
  +\mathcal{O}\left([k(T)]^{-s/d}\right) ,\ T\to \infty,\nonumber\\\label{ucfd}\end{eqnarray}
 \noindent  in the norm of the space $L^{2}(\Omega \times \mathbb{M}_{d},\mathcal{P}\otimes d\nu).$ In particular, this convergence rate also holds in probability from Chebychev inequality.
 \end{theorem}

 The proof of Theorem \ref{th1} is based on the application of Theorem 1 in \cite{Newey97} to the finite-dimensional approximation of $X$ obtained from equations (\ref{eqkshexpansion}) and (\ref{eqkshexpansionbb}) (see also equation (\ref{aaa}) in \ref{proofth1}).
\begin{remark}
\label{remmm}
 The sup--norm  convergence rate $$\mathcal{O}_{P}\left(\zeta_{k(T),M(T)}\left[\frac{\sqrt{k(T)}}{\sqrt{M(T)}}+k(T)^{-(s-d/2)/d}\right]\right)$$\noindent  can be  obtained,  in a similar way to  Theorem 1 in \cite{Newey97},  by applying Dominated Convergence Theorem. Here,  $\zeta_{k(T),M(T)}=\sup_{z\in \mathbb{M}_{d}}\|b_{k(T)}(z)\|,$ with $b_{k(T)}(z)$ being introduced in equation (\ref{eqvectsb}), and $\|\cdot\|$ denoting the Euclidean norm.
\end{remark}

\begin{remark}
\label{remconiiiapp}
Under the formulation of  (iii) in Remark \ref{condiii}, i.e., considering $M(T)=T^{\gamma },$ and $k(T)=T^{\widetilde{\alpha } },$ for $t=0,\dots,T-1,$
\begin{eqnarray}
&&\hspace*{-2cm}\int_{\Omega \times \mathbb{M}_{d}}[X_{t}(z,\xi)-\widehat{\mathbf{X}}_{t,k(T),M(T)}(z)]^{2}d\nu(z) \mathcal{P}(d\xi)\nonumber\\
&&=\mathcal{O}\left(\frac{1}{T^{\gamma -\widetilde{\alpha } }}+T^{-\widetilde{\alpha } (2s/d)}\right),\quad T\to \infty.
\label{l2cr}
\end{eqnarray}
\noindent   For $\gamma >2\widetilde{\alpha } >0,$  the rate of convergence to zero in equation  (\ref{l2cr}) is given by $T^{-\min\left\{\gamma -\widetilde{\alpha },  2\widetilde{\alpha } s/d\right\}}.$  Thus, for $t=0,\dots,T-1,$
\begin{equation}
\widehat{\mathbf{X}}_{t,k(T),M(T)}(x)\underset{L^{2}(\Omega \times \mathbb{M}_{d},\mathcal{P}\otimes d\nu)}{=}X_{t}(x)+\mathcal{O}\left(T^{-\min\left\{(\gamma -\widetilde{\alpha })/2, \widetilde{\alpha } s/d\right\}}\right), \quad T\to \infty,
\label{eqapl2}
\end{equation}
\noindent where, as before $\underset{L^{2}(\Omega \times \mathbb{M}_{d},\mathcal{P}\otimes d\nu)}{=}$ denotes the equality in the norm of the space $L^{2}(\Omega \times \mathbb{M}_{d},\mathcal{P}\otimes d\nu).$
\end{remark}

In what follows,  the notation $\widehat{\mathbf{X}}_{t,k(T),M(T)}(x)$ is simplified by considering  $\widehat{X}_{t}(x),$ where dependence on
$k(T),$ $M(T),$ and bold letter are omitted. In the subsequent development, we will work under  the formulation of  (iii) in Remark \ref{condiii}.

\section{Asymptotic analysis of the LRD testing procedure}
\label{asymptotics}
Let us first consider  the plug--in
    fDFT $\widehat{\widetilde{X}}_{\omega }^{(T)}$  of $X,$ based on \linebreak
  $\left\{ \widehat{X}_{t},\ t=0,\dots,T-1\right\},$ given by
 \begin{eqnarray}\widehat{\widetilde{X}}_{\omega }^{(T)}(x)&=&\frac{1}{\sqrt{2\pi T}}\sum_{t=0}^{T-1}\widehat{X}_{t}(x)\exp(-i\omega t),\quad x\in \mathbb{M}_{d},\quad \omega \in [-\pi,\pi].\label{fDFTbb}
\end{eqnarray}

From (\ref{eqapl2}), under $\min\left\{\frac{\gamma-\widetilde{\alpha } }{2},\frac{\widetilde{\alpha } s}{d}\right\}>1/2$
\begin{equation}
\widehat{\widetilde{X}}_{\omega }^{(T)}(x)=\widetilde{X}_{\omega }^{(T)}(x)+\mathcal{O}\left(T^{1/2-\min\left\{(\gamma -\widetilde{\alpha })/2, \widetilde{\alpha } s/d\right\}}\right),\label{a1}
\end{equation}
\noindent in the norm of the space $L^{2}(\Omega \times \mathbb{M}_{d},\mathcal{P}\otimes d\nu).$

  The plug--in periodogram operator $\widehat{\mathcal{P}}_{\omega }^{(T)}$ has kernel \begin{equation}\widehat{p}_{\omega}^{(T)}(x,y)=\frac{1}{2\pi T}\sum_{t=0}^{T-1}\sum_{s=0}^{T-1}\widehat{X}_{t}(x)\widehat{X}_{s}(y)\exp(-i\omega [t-s]),\quad x,y\in \mathbb{M}_{d}.\label{pluperop}\end{equation}

  Again, from (\ref{eqapl2}) and  under $\min\left\{\frac{\gamma-\widetilde{\alpha } }{2},\frac{\widetilde{\alpha } s}{d}\right\}>1/2,$  as $T\to \infty,$
  \begin{equation}\left\|\widehat{p}_{\omega}^{(T)}(x,y)-p_{\omega}^{(T)}(x,y)\right\|_{L^{2}(\Omega \times \mathbb{M}_{d}^{2},\mathcal{P}\otimes d\nu\otimes d\nu)}\leq h_{1}(T)=\mathcal{O}\left(T^{1/2-\min\left\{(\gamma -\widetilde{\alpha })/2, \widetilde{\alpha } s/d\right\}}\right), \label{a2}
\end{equation}
\noindent uniformly in $\omega \in [-\pi,\pi]$ under SRD.   Otherwise, constant involved in the order $\mathcal{O}\left(T^{1/2-\min\left\{(\gamma -\widetilde{\alpha })/2, \widetilde{\alpha } s/d\right\}}\right)$
  depends on $\sqrt{|\mathbb{M}_{d}|\sum_{n\in \mathbb{N}_{0}}\Gamma(n,d)f_{n}(\omega )},$ i.e.,  it depends on $\omega \in [-\pi,\pi].$

  The  plug-in weighted periodogram operator  $\widetilde{\widehat{\mathcal{F}}}_{\omega }^{(T)}$   has kernel $\widetilde{\widehat{f}}_{\omega }^{(T)}(x,y)$ satisfying, for $T\to \infty ,$
 \begin{eqnarray}
&&\widetilde{\widehat{f}}_{\omega }^{(T)}(x,y)=\left[\frac{2\pi}{T}\right]\sum_{s=0}^{T-1}  W^{(T)}\left(\omega - \frac{2\pi s}{T} \right)
 \widehat{p}_{ \frac{2\pi  s}{T} }^{(T)}(x,y)\nonumber\\
 &&=\int_{-\pi}^{\pi}\frac{1}{B_{T}}W\left(\frac{\omega -\xi }{B_{T}}\right)\widehat{p}_{ \xi }^{(T)}(x,y)d\xi  +\mathcal{O}\left(B_{T}^{-1}T^{-1}\right),\ x,y\in \mathbb{M}_{d},\label{enpb}
\end{eqnarray}
\noindent   for every  $\omega \in [-\pi,\pi],$  where the weight function $W^{(T)}$ and $\widehat{p}_{\omega}^{(T)}(x,y)$ have  been respectively introduced after equation  (\ref{enp}), and in (\ref{pluperop}).

From (\ref{a2}),
\begin{eqnarray}&&
\left\|\widetilde{\widehat{f}}_{\omega }^{(T)}-\widehat{f}_{\omega }^{(T)}\right\|_{L^{2}(\Omega \times \mathbb{M}_{d}^{2},\mathcal{P}\otimes d\nu\otimes d\nu)}\leq h_{2}(T)=\mathcal{O}\left(T^{1/2-\min\left\{(\gamma -\widetilde{\alpha})/2, \widetilde{\alpha} s/d\right\}}B_{T}^{-1/2}\right)\nonumber\\ &&\hspace*{6cm}+\mathcal{O}\left(B_{T}^{-1}T^{-1}\right),\quad  T\to \infty,\label{a3}\end{eqnarray}
\noindent  uniformly in $\omega \in [-\pi,\pi]$ under SRD.  Otherwise, constant involved in the order $\mathcal{O}\left(T^{1/2-\min\left\{(\gamma -\widetilde{\alpha})/2,\widetilde{\alpha} s/d\right\}}B_{T}^{-1/2}\right)$
  depends on $\sqrt{\sum_{n\in \mathbb{N}_{0}}\Gamma(n,d)f_{n}(\omega )},$ i.e.,  it depends on $\omega \in [-\pi,\pi].$

Our plug-in test statistics $\widehat{\mathcal{S}}_{B_{T}}$ is then formulated as
\begin{equation}\widehat{\mathcal{S}}_{B_{T}}=\sqrt{B_{T}T}\int_{[-\sqrt{B_{T}}/2, \sqrt{B_{T}}/2]}\widetilde{\widehat{\mathcal{F}}}_{\omega }^{(T)}
\frac{d\omega}{\sqrt{B_{T}}}= \int_{-\pi}^{\pi}\frac{\sqrt{B_{T}T}\widetilde{\widehat{\mathcal{F}}}_{\omega }^{(T)}\mathbb{I}_{[-\sqrt{B_{T}}/2, \sqrt{B_{T}}/2]}(\omega )}{\sqrt{B_{T}}}d\omega.\label{oslrdtintrobb}
\end{equation}

\subsection{Theoretical results under the null}
The following result provides the asymptotic probability distribution of $\widehat{\mathcal{S}}_{B_{T}}$  under $H_{0}.$
\begin{proposition}
\label{pr1GH0}
 Under conditions of Theorem \ref{th1}, and Theorem 2.2 in  \cite{Ruiz-MedinaCrujeiras24}, assume that  $\min\left\{(\gamma -\widetilde{\alpha})/2, \widetilde{\alpha} s/d\right\}>1,$ then

     \begin{equation}\widehat{\mathcal{S}_{B_{T}}}-E[\widehat{\mathcal{S}_{B_{T}}}]\to_{D} Y_{0}^{(\infty)},\quad T\to \infty,\label{lemcdrso}\end{equation}
\noindent where $\widehat{\mathcal{S}_{B_{T}}}$ has been introduced in (\ref{oslrdtintrobb}), and $Y_{0}^{(\infty)}$ is a zero--mean  Gaussian random element
in  the space  $\mathcal{S}(L^{2}(\mathbb{M}_{d},d\nu, \mathbb{R}))$ of Hilbert--Schmidt operators on $L^{2}(\mathbb{M}_{d},d\nu, \mathbb{R}),$ with autocovariance operator  $\mathcal{R}_{Y_{0}^{(\infty)}}= E\left[ Y_{0}^{(\infty)}\otimes Y_{0}^{(\infty)}\right]$ having kernel  $r_{ Y_{0}^{(\infty)}}(x_{1},y_{1}, x_{2},y_{2})$ given by\begin{eqnarray}
&&r_{Y_{0}^{(\infty)}}(x_{1},y_{1}, x_{2},y_{2})=2\pi \|W\|_{L^{2}(\mathbb{R})}^{2}\left[
f_{0}(x_{1},x_{2})f_{0}(y_{1},y_{2})\right.
\nonumber\\
&&
\left.\hspace*{1.5cm}+f_{0}(x_{1},y_{2})f_{0}(y_{1},x_{2})\right],\quad (x_{i},y_{i})\in \mathbb{M}_{d}^{2},\ i=1,2,
\label{cobk}
\end{eqnarray}
\noindent where $f_{0}(x,y)$ denotes the kernel of the spectral density operator $\mathcal{F}_{0}$ at frequency $\omega =0$ (see equation (\ref{sdo2})).
  \end{proposition}
The proof follows straightforward from equations (\ref{a1}), (\ref{a2}) and (\ref{a3}), under conditions assumed in  Theorem \ref{th1}, and Theorem 2.2 in  \cite{Ruiz-MedinaCrujeiras24} (see \ref{pr1}).
\subsection{Second and fourth order spectral  bias asymptotics  under $H_{1}$}
The following lemmas and results provide  the $L^{2}$ convergence rates to zero of the bias of the integrated   plug--in empirical second and fourth order
spectral cumulant operators under the conditions of Theorem \ref{th1}.

\begin{lemma}
\label{lem1}
 Under  $H_{1},$ assuming  the conditions of Theorem \ref{th1} hold. \linebreak If $\min\left\{(\gamma -\widetilde{\alpha})/2, \widetilde{\alpha} s/d\right\}>1/2,$ then
\begin{eqnarray}&&
\left\|\int_{-\pi}^{\pi}E_{H_{1}}\left[\widehat{\mathcal{P}}_{ \omega }^{(T)} \right]d\omega -\int_{-\pi}^{\pi}\mathcal{F}_{\omega }d\omega\right\|_{\mathcal{S}(L^{2}(\mathbb{M}_{d},d\nu, \mathbb{C}))}\nonumber\\
&&\hspace*{1cm}\leq h_{3}(T)=
 \mathcal{O}(T^{-1})+\mathcal{O}\left(T^{1/2-\min\left\{(\gamma -\widetilde{\alpha})/2, \widetilde{\alpha} s/d\right\}}\right), \quad T\to \infty,\nonumber\\
 \label{orderofconvmpo}
   \end{eqnarray}

 \noindent   where $E_{H_{1}}$ denotes expectation under the  alternative $H_{1},$ with
$\mathcal{S}(L^{2}(\mathbb{M}_{d},d\nu, \mathbb{C}))\equiv L^{2}(\mathbb{M}_{d}^{2}, d\nu\otimes d\nu,\mathbb{C}).$
\end{lemma}
The proof of this lemma follows from triangle inequality, applying  Lemma 3.1 in \cite{Ruiz-MedinaCrujeiras24}, equation (\ref{a2}), and Jensen's inequality, allowing the application of  Dominated Convergence Theorem (see \ref{lem1bb}).

\begin{corollary}\label{lem3}
Under $H_{1},$ assume that $B_{T}=T^{-\beta },$ $\beta \in (0,1),$ and  the conditions of Theorem \ref{th1} hold.
If  $\min\left\{(\gamma -\widetilde{\alpha})/2, \widetilde{\alpha} s/d\right\}>(1+\beta )/2,$ then
\begin{eqnarray}&&\hspace*{-1cm}\left\|\int_{-\pi}^{\pi}E_{H_{1}}\left[\widetilde{\widehat{\mathcal{F}}}^{(T)}_{\omega }\right] -\int_{\mathbb{R}}W(\xi)\mathcal{F}_{\omega -\xi B_{T}}d\xi d\omega \right\|_{\mathcal{S}(L^{2}(\mathbb{M}_{d},d\nu, \mathbb{C}))}\nonumber\\ &&\leq h_{4}(T)=\mathcal{O}(B_{T}^{-1}T^{-1})+\mathcal{O}(T^{-1})+\mathcal{O}\left(T^{1/2-\min\left\{(\gamma -\widetilde{\alpha} )/2, \widetilde{\alpha} s/d\right\}}B_{T}^{-1/2}\right) \nonumber\\ &&\hspace*{-0.5cm}=\mathcal{O}(B_{T}^{-1}T^{-1})+\mathcal{O}(T^{-1})+\mathcal{O}\left(T^{(1+\beta )/2-\min\left\{(\gamma -\widetilde{\alpha} )/2, \widetilde{\alpha} s/d\right\}}\right),\quad T\to \infty.\label{eqapphh}\end{eqnarray}
\end{corollary}

The proof of Corollary \ref{lem3} follows from triangle inequality, applying  Corollary 3.2 in \cite{Ruiz-MedinaCrujeiras24}, Jensen's inequality, and equation (\ref{a3}), allowing the application of Dominated Convergence Theorem (see \ref{apcor3}).

\begin{lemma}
\label{lem4cs}
Under  $H_{1},$ assume that  for $i\neq j\neq k,$ with $i,j,k\in \{1,2,3,4\},$
\begin{equation}\int_{[-\pi,\pi]^{3}}
E\left\|\widetilde{X}_{\omega_{i}}\otimes \widetilde{X}_{\omega_{j}}\otimes \widetilde{X}_{\omega_{k}}\right\|^{2}_{L^{2}(\mathbb{M}_{d}^{3},\otimes_{i=1}^{3}d\nu_{i})}d\omega_{i}d\omega_{j}d\omega_{k}<\infty.
\label{eqlw33}
\end{equation}
Considering  the conditions in Lemma 3.3  in \cite{Ruiz-MedinaCrujeiras24},   and in Theorem \ref{th1} hold, if $\min\{(\gamma -\widetilde{\alpha} )/2,\widetilde{\alpha} s/d\}>3/2,$ then
\begin{eqnarray}&&
\left\|\int_{[-\pi ,\pi ]^{3}}
T\mbox{cum}\left(\widehat{\widetilde{X}}_{\omega_{1}}^{(T)}(\tau_{1}), \widehat{\widetilde{X}}_{\omega_{2}}^{(T)}(\tau_{2}),
\widehat{\widetilde{X}}_{\omega_{3}}^{(T)}(\tau_{3}), \widehat{\widetilde{X}}_{\omega_{4}}^{(T)}(\tau_{4})\right)\right.\nonumber\\
&&\left.\hspace*{1cm}-2\pi \int_{[-\pi ,\pi ]^{3}}\mathcal{F}_{\omega_{1},\omega_{2},\omega_{3}}(\tau_{1},\tau_{2}, \tau_{3}, \tau_{4})
d\omega_{1}d\omega_{2}d\omega_{3}\right\|_{\mathcal{S}\left(L^{2}\left(\mathbb{M}_{d}^{2},\otimes_{i=1}^{2}\nu(dx_{i}), \mathbb{C}\right)\right)}\nonumber\\
&&\leq h_{5}(T)=\mathcal{O}(T^{-1})+\mathcal{O}\left(T^{3/2-\min\left\{(\gamma -\widetilde{\alpha} )/2, \widetilde{\alpha} s/d\right\}}\right),\quad T\to \infty ,
\label{eqapp}\end{eqnarray}
\noindent
 uniformly in $\omega_{4}\in [-\pi,\pi],$ where, for $\omega_{i}\in [-\pi,\pi],$ $i=1,2,3,4,$ $\widehat{\widetilde{X}}_{\omega_{i}}^{(T)}$ denotes, as before, the plug--in
    fDFT introduced in equation (\ref{fDFTbb}),  and
\begin{eqnarray}&&\hspace*{-0.35cm} \mathcal{F}_{\omega_{1}, \omega_{2}, \omega_{3}}\underset{\mathcal{S}\left(L^{2}\left(\mathbb{M}_{d}^{2},\otimes_{i=1}^{2}\nu(dx_{i}), \mathbb{C}\right)\right)}{=}\frac{1}{(2\pi)^{3}}\sum_{t_{1},t_{2},t_{3}= -\infty}^{\infty}
\exp\left(\sum_{j=1}^{3}\omega_{j}t_{j}\right)\nonumber\\
&&\hspace*{7cm}\times
\mbox{cum}(X_{t_{1}},X_{t_{2}},X_{t_{3}},X_{0})\nonumber\\
\label{csdoo4}
\end{eqnarray}

\noindent
 is the cumulant spectral density operator of order $4$ of $X.$  Here, \linebreak $\underset{\mathcal{S}\left(L^{2}\left(\mathbb{M}_{d}^{2},\otimes_{i=1}^{2}\nu(dx_{i}), \mathbb{C}\right)\right)}{=}$
means the identity in the norm of the  space \linebreak $\mathcal{S}\left(L^{2}\left(\mathbb{M}_{d}^{2},\otimes_{i=1}^{2}\nu(dx_{i}), \mathbb{C}\right)\right)\equiv L^{2}(\mathbb{M}^{4}_{d}, \otimes_{i=1}^{4} d\nu (x_{i}),\mathbb{C}).$

\end{lemma}
The proof of Lemma \ref{lem4cs} follows from Lemma 3.3 in \cite{Ruiz-MedinaCrujeiras24}, and equation (\ref{a1}), applying fourth order cumulant formula  in terms of moments, Jensen's inequality and Dominated Convergence Theorem, in a similar way to the proof of the previous lemmas.

\subsection{Consistency}
\label{sconsistency}
Consistency of the test is based on the a.s. divergence of the plug--in test statistics $\widehat{\mathcal{S}}_{B_{T}}$ in the Hilbert-Schmidt operator norm.  The next two results play a crucial role in the derivation of consistency.
\begin{proposition}
\label{cor1}
Under $H_{1},$  assume that the conditions of Theorem \ref{th1} hold. Then,
\begin{eqnarray}&&\hspace*{-0.7cm}\left\|E_{H_{1}}\left[\frac{\widehat{\mathcal{S}}_{B_{T}}}{\sqrt{TB_{T}}}\right]\right\|_{\mathcal{S}(L^{2}(\mathbb{M}_{d},d\nu, \mathbb{C}))}
=
\left\|\int_{[-\sqrt{B_{T}}/2, \sqrt{B_{T}}/2]}E_{H_{1}}[\widetilde{\widehat{\mathcal{F}}}_{\omega }^{(T)}]\frac{d\omega}{\sqrt{B_{T}}} \right\|_{\mathcal{S}(L^{2}(\mathbb{M}_{d},d\nu, \mathbb{C}))}\nonumber\\
&&\hspace*{4cm}
\geq g(T)=\mathcal{O}(B_{T}^{-l_{\alpha }-1/2}),\quad  T\to \infty.\label{eqorexh1hh}\end{eqnarray}
\end{proposition}

\noindent The proof of this result is obtained applying the properties of F\'ejer kernel introduced in equation  (\ref{fejkernel}), of the weighted function  $W^{(T)}$ in equation (\ref{enp}), and the boundedness of eigenvalues of operator $\mathcal{A}$ in equation (\ref{lerdop}) (see \ref{pproposition2}).

\begin{theorem}
\label{consisimse} Under $H_{1},$ considering  the conditions of Theorem \ref{th1},  Corollary \ref{lem3}, and  Theorem 4.2 in  \cite{Ruiz-MedinaCrujeiras24} hold,
\begin{eqnarray}
&&\hspace*{-0.7cm}\int_{-\pi}^{\pi}E_{H_{1}}\left\|\widetilde{\widehat{\mathcal{F}}}_{\omega }^{(T)}-E_{H_{1}}[\widetilde{\widehat{\mathcal{F}}}_{\omega }^{(T)}]\right\|_{L^{2}(\mathbb{M}_{d}^{2},d\nu\otimes d\nu, \mathbb{C})}^{2}d\omega \leq h_{6}(T)\nonumber\\ &&\hspace*{0.7cm}=\mathcal{O}(B_{T}^{-1}T^{-1})+\mathcal{O}\left(T^{(1+\beta)-2\min\left\{(\gamma -\widetilde{\alpha })/2  , \widetilde{\alpha } s/d\right\}}\right),\ T\to \infty.
\label{ca2}
\end{eqnarray}

\end{theorem}

\noindent he proof of Theorem \ref{consisimse} is obtained from  applying the  triangle, Cauchy--Schwartz,  and Jensen's inequalities, considering   equation (\ref{a3}),    Dominated Convergence Theorem,   and Theorem 4.2 in  \cite{Ruiz-MedinaCrujeiras24} (see \ref{consist}).

  The   weak consistency of the  integrated plug--in  weighted periodogram operator   under $H_{1}$ is obtained in the next result, in the norm of    the space $L^{2}(\mathbb{M}_{d}^{2},d\nu\otimes d\nu, \mathbb{C})$. This result follows from  Corollary \ref{lem3} and Theorem \ref{consisimse}, applying triangle inequality,  in a similar way to Corollary 4.4 in \cite{Ruiz-MedinaCrujeiras24}.
\begin{corollary}
\label{cor22}
Under the conditions of Theorem \ref{consisimse}, as $T\to \infty,$
 \begin{eqnarray}
&&\left\|\int_{-\pi}^{\pi}E_{H_{1}}\left[\widetilde{\widehat{\mathcal{F}}}_{\omega }^{(T)}-\int_{-\pi}^{\pi}W(\xi)\mathcal{F}_{\omega -B_{T}\xi }d\xi \right] d\omega
\right\|_{\mathcal{S}(L^{2}(\mathbb{M}_{d},d\nu ,\mathbb{C}))}\nonumber\\ &&\leq h_{7}(T)=\mathcal{O}(T^{-1/2}B_{T}^{-1/2})+\mathcal{O}\left(T^{(1+\beta)/2-\min\left\{(\gamma -\widetilde{\alpha } )/2  , \widetilde{\alpha } s/d\right\}}\right).
\nonumber\end{eqnarray}
\end{corollary}
The following theorem provides  consistency of the  test under   strong correlated in time functional data, that are  (irregularly) discretely observed, and  affected by additive spatiotemporal white noise.
\begin{theorem}
\label{th2}
Assume that  $H_{1}$ holds  with $l_{\alpha }> 1/4.$ Under  the conditions of Theorem \ref{consisimse}, with  $B_{T}=T^{-\beta },$  $\beta \in (0,1),$
 if the parameters $s, \widetilde{\alpha }, \gamma $ and $\beta $
are  such that
\begin{equation}\min\left\{\gamma-\widetilde{\alpha }, 2s\widetilde{\alpha } / d\right\}>  2-\beta \left[2l_{\alpha }-\frac{1}{2}\right],\label{cf}\end{equation}
\noindent then,
$$
\left\|\widehat{\mathcal{S}}_{B_{T}}\right\|_{\mathcal{S}(L^{2}(\mathbb{M}_{d},d\nu, \mathbb{C}))}\to_{\mbox{a.s}} \infty,$$
\noindent where $\to_{\mbox{a.s.}} \infty$ denotes a.s. divergence.
\end{theorem}

\noindent The proof of Theorem \ref{th2}  follows from
  Proposition   \ref{cor1}  and Theorem \ref{consisimse}, applying a  similar  methodological approach   to the one considered in the proof of Theorem 4.5 in \cite{Ruiz-MedinaCrujeiras24} (see \ref{maresult}).

\section{Simulation study. Practical illustrations}
\label{simulationstudy}
 Along this section we consider  $B_{T}=T^{-\beta },$ $\beta =1/4.$ For a deeper analysis of $\beta $ parameter values (see
Section 6 in \cite{Ruiz-MedinaCrujeiras24}).
The consistency of the test is illustrated showing the a.s. divergence of the Hilbert--Schmidt operator norm of the plug--in test statistics $\widehat{S}_{B_{T}},$ under  three LRD operator models corresponding to  Examples 1--3, which display different locations of the dominant subspace contained into an eigenspace of the Laplace Beltrami operator. The empirical test size and power are also computed for different functional sample sizes $T.$
\subsection{Consistency}
\label{simuconsistency}
As commented, the numerical results in this section  illustrate the a.s. divergence of the Hilbert--Schmidt operator norm of the plug--in test statistics $\widehat{S}_{B_{T}},$ under uniform spherical   sparse spatial sampling design.  Data
in Examples 1--3 have been generated, conditionally to the uniform spatial sampling design,  from multifractionally integrated SPHARMA(p,q)  models  (see Sections 5.2.1--5.2.3 in \cite{Ruiz-MedinaCrujeiras24} for more details). These data set is  contaminated with additive observation noise $\varepsilon $  having variance  $\sigma^{2}=1/8.$
In condition (ii) we have considered $s=3.$

\begin{figure}[!h]
\begin{center}
\includegraphics[width=2.5cm,height=2.5cm]{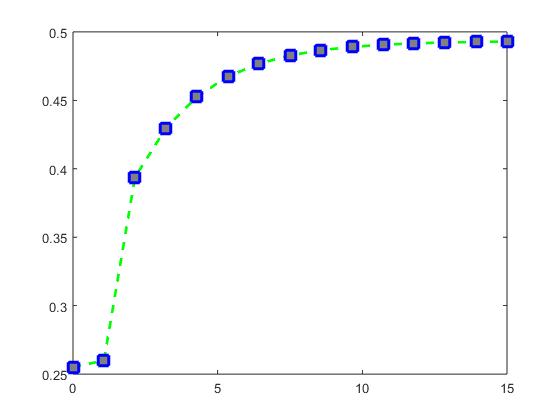}
\includegraphics[width=2.5cm,height=2.5cm]{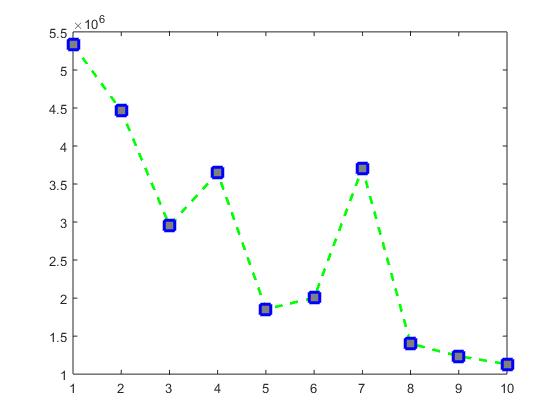}
\includegraphics[width=2.5cm,height=2.5cm]{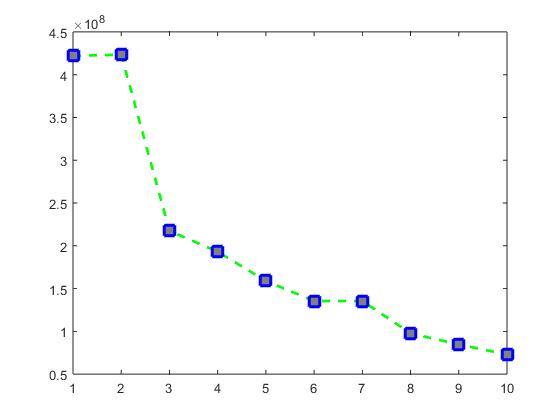}
\includegraphics[width=2.5cm,height=2.5cm]{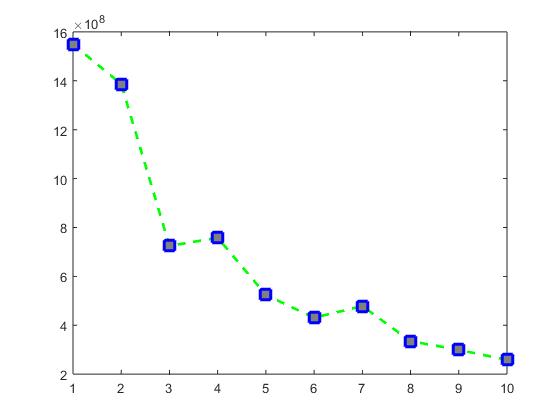}
\includegraphics[width=2.5cm,height=2.5cm]{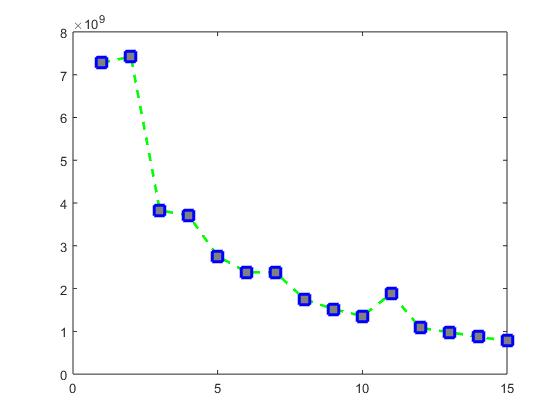}
\end{center}
\caption{\emph{Example 1.} Eigenvalues $\alpha (n,j),$ $j=1,\dots,\Gamma (n,2),$ $n=1,\dots, N(k(T)),$  of  LRD operator $\mathcal{A},$  for $N(k(T))=5,$ and $k(T)= 15,$ with  $L_{\alpha }= 0.4929,$ and $l_{\alpha }=0.2550$  (left--hand side).  Sample  projections of  $\widehat{\mathcal{S}}_{B_{T}},$ $B_{T}= T^{-1/4},$
into the eigenfunctions generating the tensor product  eigenspaces   $\mathcal{H}_{n}\otimes \mathcal{H}_{n},$ $n=1,\dots,N(k(T))$    (four plots at the right--hand--side)}
\label{CONS1}
\end{figure}

\begin{figure}[!h]
\begin{center}
\includegraphics[width=2.5cm,height=2.5cm]{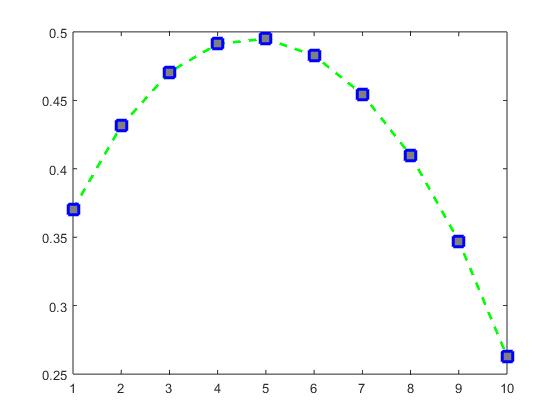}
\includegraphics[width=2.5cm,height=2.5cm]{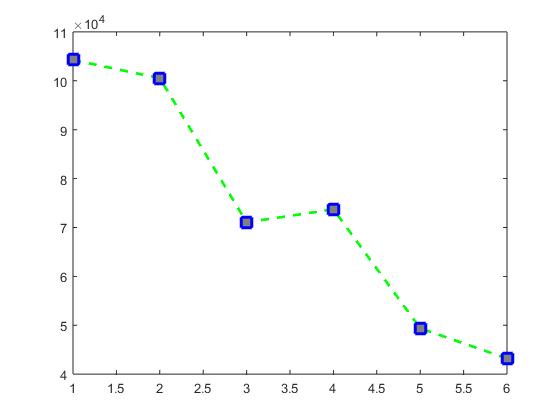}
\includegraphics[width=2.5cm,height=2.5cm]{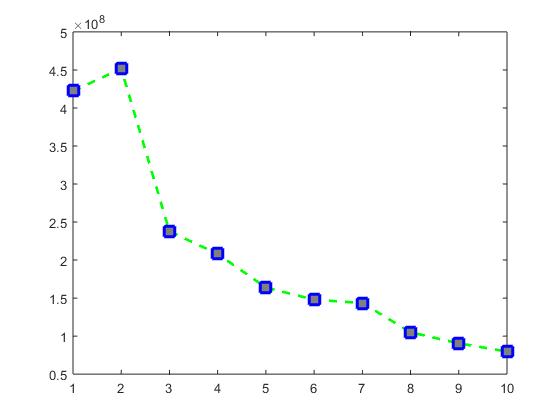}
\includegraphics[width=2.5cm,height=2.5cm]{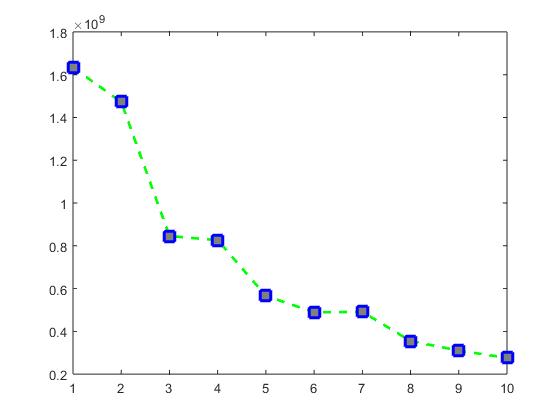}
\includegraphics[width=2.5cm,height=2.5cm]{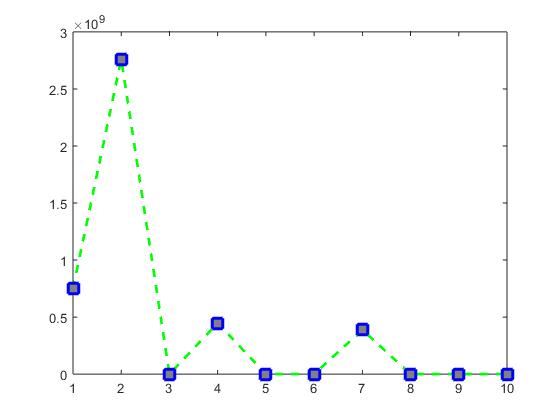}
\end{center}
\caption{\emph{Example 2.} Eigenvalues $\alpha (n,j),$ $j=1,\dots,\Gamma(n,2),$ $n=1,\dots, N(k(T)),$    of  LRD operator $\mathcal{A},$
 for $N(k(T))=4,$ and $k(T)= 10,$
with $L_{\alpha }= 0.4950,$ and $l_{\alpha }=0.2629$    (left--hand side).  Sample  projections of   $\widehat{\mathcal{S}}_{B_{T}},$ $B_{T}= T^{-1/4},$
into the eigenfunctions generating the tensor product eigenspaces   $\mathcal{H}_{n}\otimes \mathcal{H}_{n},$ $n=1,\dots,N(k(T)),$  $T=1000,30000,50000,100000$ (four plots at the right--hand--side)}
\label{CONSE1EX2}
\end{figure}

\begin{figure}[!h]
\begin{center}
\includegraphics[width=2.5cm,height=2.5cm]{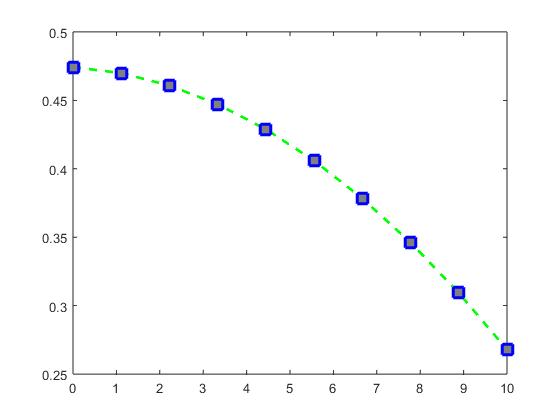}
\includegraphics[width=2.5cm,height=2.5cm]{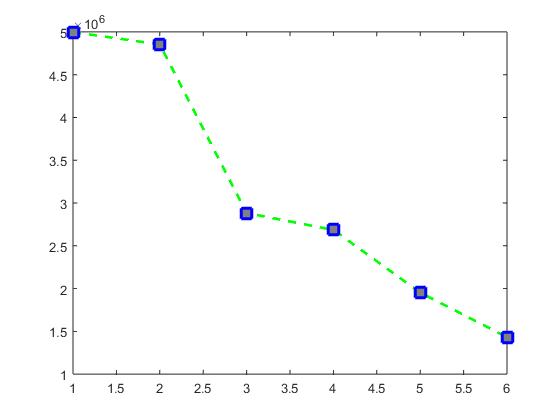}
\includegraphics[width=2.5cm,height=2.5cm]{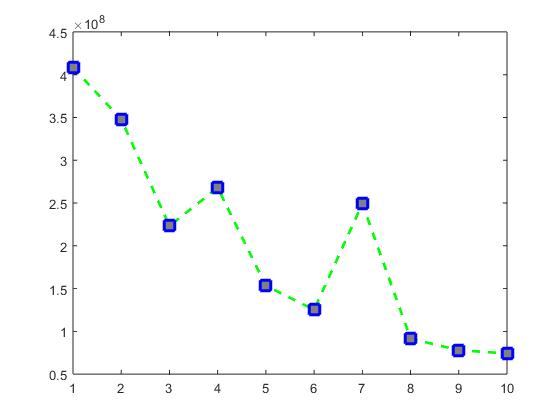}
\includegraphics[width=2.5cm,height=2.5cm]{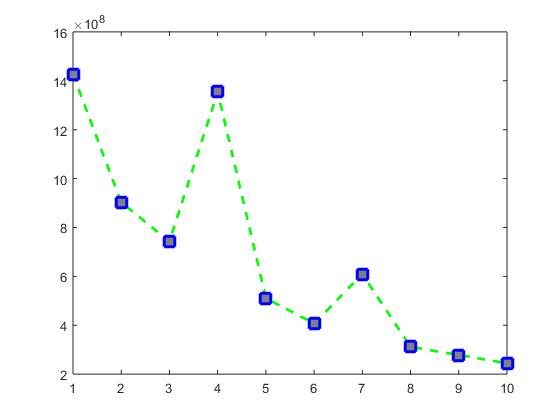}
\includegraphics[width=2.5cm,height=2.5cm]{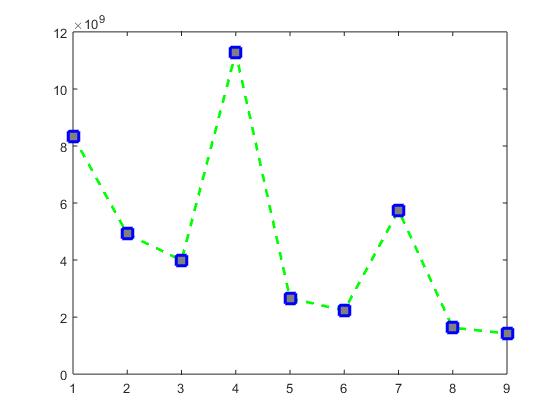}
\end{center}
\caption{\emph{Example 3.} Eigenvalues $\alpha (n,j),$ $j=1,\dots, \Gamma (n,2),$ $n=1,\dots, N(k(T)),$   of  LRD operator $\mathcal{A},$  for $N(k(T))=4,$ and $k(T)= 10,$ with $L_{\alpha }=   0.4743,$ and $l_{\alpha }= 0.2678$    (left--hand side).  Sample projections of  $\widehat{\mathcal{S}}_{B_{T}},$ $B_{T}= T^{-1/4},$
into  the eigenfunctions generating the tensor product eigenspaces $\mathcal{H}_{n}\otimes \mathcal{H}_{n},$ $n=1,\dots,N(k(T)),$ $T=1000,30000,50000,100000$  (four plots at the right--hand--side)}
\label{CONSE1EX3}
\end{figure}

Under $H_{1},$ we consider the  frequency--varying  eigenvalues $\left\{ f_{n}(\omega ),\ n\in \mathbb{N}_{0}\right\}$ satisfying
\begin{eqnarray}
\hspace*{-0.5cm} f_{n}(\omega ) &=&\frac{\lambda_{n}(\mathcal{R}^{\eta}_{0})}{2\pi}\left|\frac{\Psi_{q,n}(\exp(-i\omega ))}{\Phi_{p,n}(\exp(-i\omega ))}\right|^{2}\left|1-\exp\left(-i\omega \right)\right|^{-\alpha (n,j)},\ n\in \mathbb{N}_{0},
\label{eqsim2bb}
\end{eqnarray}
\noindent   for each $\omega \in [-\pi,\pi],$ where $\left\{\alpha (n,j),\ j=1,\dots,\Gamma (n, 2), \ n\in \mathbb{N}_{0}\right\}$ define the eigenvalues of the LRD operator $\mathcal{A}$  in  Examples 1--3, respectively plotted  at the left--hand side of  Figures \ref{CONS1}--\ref{CONSE1EX3} below (see also Sections 5.2.1--5.2.3 in \cite{Ruiz-MedinaCrujeiras24}, for more details on the generated  multifractionally integrated SPHARMA(1,1) functional time series models). Data spherical sparsity is controlled by the shape  parameter value $\gamma = 0.3077,$ under $M(T)=\left[[T^{-\gamma}]^{-}\right]^{2},$  with $[\cdot ]^{-}$ denoting the integer part function. Local smoothness parameter value  $s=3,$  and the specified  values below of $T,$ $k(T)$ and $M(T)$ lead to our above choice of  value
  $\sigma^{2}=1/8,$ for a suitable noise to signal ratio.

  In Example 1, we consider the dominant subspace is included into the $N(k(T))$--th eigenspace of the Laplace Beltrami operator,  being  generated by the  eigenfunctions  associated with the $N(k(T))$--th  eigenvalue of the Laplace Beltrami operator  (see left--hand side of Figure \ref{CONS1}, with  $L_{\alpha }=  0.4929,$  $l_{\alpha }=  0.2550,$ and   $\alpha (n,j)=L_{\alpha }= 0.4929,$ $j=1\dots,\Gamma(n,2),$ $n\geq 15$). As before, $N(k(T))$ is such that  $\sum_{l=0}^{N(k(T))}\Gamma (l,2)=k(T).$ Considering    $B_{T}= T^{-1/4},$ the sample values   of the projections of  $\widehat{\mathcal{S}}_{B_{T}}$
  into the $k(T)$ eigenfunctions, generating the tensor product  eigenspaces   $\mathcal{H}_{n}\otimes \mathcal{H}_{n},$ $n=1,\dots,N(k(T)),$   are plotted at the right--hand side of Figure \ref{CONS1}  (see also Table \ref{T1} where the corresponding  truncated Hilbert--Schmidt operator norms of $\widehat{S}_{B_{T}}$ are displayed).

In the second example,  since, as given at the left--hand side of  Figure \ref{CONSE1EX2},  the dominant subspace is generated by the fifth eigenfunction associated with the third eigenvalue of the Laplace Beltrami operator (i.e., this subspace is contained into the third eigenspace of the Laplace Beltrami operator), lower sieve basis dimension is required to include  the eigenfunction generating this dominant subspace. Hence, we consider as maximum sieve basis dimension $k(T)=[\log(T)-1]^{-}=10,$ for $T=100000,$ being $N(k(T))=4.$ Figure \ref{CONSE1EX2} displays, at the right--hand side, the sample projections of  $\widehat{\mathcal{S}}_{B_{T}},$ $B_{T}= T^{-1/4},$
into the eigenfunctions generating the tensor product eigenspaces   $\mathcal{H}_{n}\otimes \mathcal{H}_{n},$ $n=1,\dots,N(k(T)),$ $T=1000,30000,50000,100000$ (see also Table \ref{T1}, where the truncated norms are displayed).

Finally, in the third example, we consider that the dominant subspace is generated by the first eigenfunction associated with the first eigenvalue of the Laplace Beltrami operator (see Figure \ref{CONSE1EX3} at the left--hand side). We then consider, as in the previous example, as maximum sieve basis dimension $k(T)=[\log(T)-1]^{-}=10,$ for $T=100000,$ being $N(k(T))=4.$   The sample  projections of  $\widehat{\mathcal{S}}_{B_{T}},$ $B_{T}= T^{-1/4},$
into the eigenfunctions generating the  tensor product  eigenspaces $\mathcal{H}_{n}\otimes \mathcal{H}_{n},$ $n=1,\dots,N(k(T)),$ $T=1000,30000,50000,100000,$ are displayed in the four plots at the right--hand--side of Figure \ref{CONSE1EX3}  (see also Table \ref{T1}).

\begin{table*}
\caption{\textbf{\emph{Hilbert-Schmidt operator norm of projected $\widehat{\mathcal{S}}_{B_{T}},$ $\beta =1/4$} }}
\label{T1}
\medskip
\begin{center}
\begin{tabular}{@{}lrrrrc@{}}
\hline
 {\bf T}& \multicolumn{1}{c}{{\bf Example 1}} &\multicolumn{1}{c}{{\bf Example 2}}    &\multicolumn{1}{c}{{\bf Example 3}} \\
  \hline
 1000 &   1.8669e+05 &      1.8911e+05    &      1.1680e+06\\
5000 &  9.8156e+06  &       1.0445e+07&    8.3620e+06\\
10000 &     6.8435e+07 &    7.0899e+07  &   5.0493e+07 \\
30000&  7.2519e+08&     7.5960e+08 &   7.2832e+08\\
50000&  2.5255e+09 &  2.7079e+09 &  2.5032e+09 \\
100000 & 1.3027e+10 &  2.9231e+09 &    1.6924e+10 \\
  \hline
\end{tabular}
\end{center}
\end{table*}

\subsection{Empirical size and power analysis}
\label{secspowt}

The results displayed in this section are based on $R=500$ repetitions, since according to the  parameter values of $T,$ $k(T),$ and $M(T)$ analyzed  ($M(T)=T^{2\gamma },$ $\gamma \in (0,1),$    $k(T)\simeq \mathcal{O}\left(\log(T)\right),$ $T\to \infty$), larger values of  $R$ do not   improve the results obtained.
  The empirical  test sizes and powers  are computed considering  six random projections (see \cite{Ruiz-MedinaCrujeiras24} for more details on the conditions required to apply  the random projection methodology in this context).

Under $H_{0}$, to compute the empirical test size, we consider the elements of the spectral density operator family $\{ \mathcal{F}_{\omega },\ \omega \in [-\pi,\pi]\}$ have pure point spectra given by (see \cite{Ruiz-MedinaCrujeiras24})
\begin{eqnarray}
f_{n}(\omega ) &=&\frac{\lambda_{n}(\mathcal{R}^{\eta}_{0})}{2\pi}\left|\frac{\Psi_{q,n}(\exp(-i\omega ))}{\Phi_{p,n}(\exp(-i\omega ))}\right|^{2},\ n\in \mathbb{N}_{0},\ \omega \in [-\pi,\pi],\nonumber\\
\label{eqsim1}
\end{eqnarray}
\noindent where $\left\{\lambda_{n}(\mathcal{R}^{\eta}_{0}),\ n\in \mathbb{N}_{0}\right\}$ is the system of eigenvalues of the trace  autocovariance operator $\mathcal{R}^{\eta}_{0}$ of the innovation process $\eta=\{\eta_{t},\ t\in \mathbb{Z}\},$ with respect to the system of eigenfunctions of the Laplace--Beltrami operator.  Process   $\eta $ is assumed to be strong--white noise in $L^{2}(\mathbb{S}_{d},d\nu, \mathbb{R}).$ That is, $\eta $ is assumed to be a
sequence of independent and identically distributed $L^{2}(\mathbb{S}_{d},d\nu, \mathbb{R})$--valued random variables
such that $E[\eta_{t}]=0,$ and $E[\eta_{t}\otimes
\eta_{s}]=\delta_{t,s}\mathcal{R}^{\eta}_{0},$  and $\delta_{t,s}=0,$
for $t\neq s,$ and $\delta_{t,s}=1,$ for $t=s.$
For $n\in \mathbb{N}_{0},$
  $\Phi_{p,n}(z)=1-\sum_{j=1}^{p}\lambda_{n}(\varphi_{j})z^{j}$ and $\Psi_{q,n}(z)=\sum_{j=1}^{q}\lambda_{n}(\psi_{j})z^{j},$
 with $\left\{\lambda_{n}(\varphi_{j}),\ n\in \mathbb{N}_{0}\right\}$ and  $\left\{\lambda_{n}(\psi_{l}),\ n\in \mathbb{N}_{0}\right\}$
     denoting the sequences of eigenvalues, with respect to the system of eigenfunctions of the Laplace--Beltrami operator,  of the  self--adjoint invariant integral operators $\varphi_{j}$   and $\psi_{l},$ for $j=1,\dots,p,$ and $l=1,\dots,q,$ respectively.   These operators  satisfy the following equations:
\begin{eqnarray}
\Phi_{p}(B)=1-\sum_{j=1}^{p}\varphi_{j}B^{j},\quad \Psi_{q}(B)=\sum_{j=1}^{q}\psi_{j}B^{j},\nonumber 
\label{rf}
\end{eqnarray}
\noindent where   $B$ is the difference operator  introduced in  \cite{RuizMedina2022}. Thus, $\Phi_{p}$ and  $\Psi_{q}$ are the so--called autoregressive and moving averages operators, respectively.
 Also, for each $n\in \mathbb{N}_{0},$  $\Phi_{p,n}(z)=1-\sum_{j=1}^{p}\lambda_{n}(\varphi_{j})z^{j}$ and $\Psi_{q,n}(z)=\sum_{j=1}^{q}\lambda_{n}(\psi_{j})z^{j}$ have not common
roots, and  their roots are outside of the unit circle (see also Corollary 6.17 in \cite{Beran17}). Thus, $X$ satisfies   an SPHARMA(p,q) equation (see, e.g.,  \cite{CaponeraMarinucci}).

Table \ref{T3bb} displays the empirical test sizes obtained  for  the functional sample sizes $T=50,100,500,1000,$ based on $R=500$ repetitions,  under
\linebreak SPHARMA(1,1) process scenario, i.e.,  we consider   $p=q=1,$ with $H=L^{2}(\mathbb{S}_{2},d\nu, \mathbb{R}),$  and $\lambda_{n}(\varphi_{1})=0.7\left(\frac{n+1}{n}\right)^{-3/2}$ and
$\lambda_{n}(\psi_{1})= (0.4)\left(\frac{n+1}{n}\right)^{-5/1.95},$ $n\in \mathbb{N}_{0}$   (see also Section 5.1 in \cite{Ruiz-MedinaCrujeiras24}, for more details on this generated model).
For each functional sample size $T=50,100,500, 1000$,  the empirical test size $\widehat{\alpha }$ is computed   for different spatial sampling frequencies $M(T)=T^{2\gamma},$ with $\gamma =0.45, 0.3704, 0.35, 0.25.$ As commented,   $k(T)\simeq [\log(T)]^{-},$  allowing the choice $\sigma^{2}=1/2$ for the  variance of the additive spatiotemporal observation white noise $\varepsilon, $ keeping in mind the displayed  $M(T)$ values.

\begin{table}[!h]
\caption{\textbf{\emph{SPHARMA(1,1). Empirical size.  $M(T)=T^{2\gamma},$ $\beta =1/4,$   $\mathbf{k}_{n,j,h,l},$ $n=h=1,2,3$ $\alpha =0.05$} }}
\begin{center}
\begin{tabular}{|l|llllll|}
\hline
   {\bf R}$=$ 500 & & & & &&   \\ \hline
   {\bf T}$=$50,  $\gamma=0.4500 $ &  0.0320  &   0.0560 &   0.0480  &  0.0360  &  0.0320   & 0.0360\\
     {\bf T}$=$50,  $\gamma=0.3704$  &  0.0280   & 0.0440 &   0.0480   & 0.0280   & 0.0400 &   0.0440\\
     {\bf T}$=$50,  $\gamma=0.3500$ &  0.0240  &   0.0480  &   0.0640  &   0.0320  &   0.0480 &    0.0320\\
  {\bf T}$=$50,  $\gamma=0.2500$ & 0.0120  &   0.0560 &    0.0640   &  0.0160   &  0.0560   &  0.0320\\
   \hline
 {\bf T}$=$100, $\gamma=0.4500$ & 0.0560 &   0.0320 &   0.0520 &   0.0520 &   0.0280 &   0.0440\\
    {\bf T}$=$100, $\gamma=0.3704$ & 0.0520  &  0.0320  &  0.0440 &   0.0400  &  0.0600  &  0.0520\\
      {\bf T}$=$100, $\gamma=0.3500$ &   0.0520&    0.0320  &  0.0440   & 0.0400 &   0.0600  &  0.0520\\
{\bf T}$=$100, $\gamma=0.2500$ &  0.0480 &   0.0240 &   0.0280 &   0.0240  &  0.0520  &  0.0600\\
   \hline
    {\bf T}$=$500, $\gamma=0.4500$ &   0.0520  &  0.0560 &   0.0640 &   0.0520 &   0.0560 &   0.0680
\\
{\bf T}$=$500, $\gamma=0.3704$ & 0.0240 &   0.0280 &   0.0440  & 0.0520 &   0.0440 &   0.0480\\
     {\bf T}$=$500, $\gamma=0.3500$ &     0.0520 &   0.0320   & 0.0440&    0.0400  &  0.0600  &  0.0520\\
      {\bf T}$=$500, $\gamma=0.2500$ & 0.0600  &  0.0440 &   0.0640  &  0.0520  &  0.0400  &  0.0720\\
   \hline
     {\bf T}$=$1000, $\gamma=0.4500$ &  0.0280  &  0.0440 &   0.0560  &  0.0360 &   0.0400  &  0.0520\\
      {\bf T}$=$1000, $\gamma=0.3704$&  0.0360  &  0.0440  &  0.0240  &  0.0400  &  0.0440  &  0.0400\\
  {\bf T}$=$1000, $\gamma=0.3500$ &  0.0360 &   0.0400 &   0.0280  &  0.0440  &  0.0360  &  0.0520\\
    {\bf T}$=$1000, $\gamma=0.2500$ &   0.0360   &  0.0240  &   0.0360  &   0.0560  &   0.0360   &  0.0480\\
    \hline
 \end{tabular}\label{T3bb}
\end{center}
\end{table}

In the empirical power analysis,  we have considered the three LRD scenarios provided by Examples 1--3 in Section \ref{simuconsistency}.
As before, our   irregularly discretely observed spherical  functional data are contaminated by   additive spatiotemporal white noise with variance  $\sigma^{2}=1/2.$  Tables \ref{T4}--\ref{T6} provide, for the three examples analyzed in Section \ref{simuconsistency}, the computed  empirical powers for the functional samples sizes  $T=50, 100, 500,1000,$ based on   $R=500.$ For each functional sample size, we analyze the  effect of the uniform spherical sampling frequency, by considering different values of  the  parameter $\gamma $ reflecting sparsity of the spatial discrete observations of our functional data set (according to  $M(T)=T^{2\gamma},$ $\gamma \in (0,1)$).

 \begin{table}[!h]
\caption{\textbf{\emph{Example 1. Empirical power $M(T)=T^{2\gamma},$ $\beta =1/4,$   $\mathbf{k}_{n,j,h,l},$ $n=h=1,2,3,$ $\alpha =0.05$} }}
\begin{center}
\begin{tabular}{|l|llllll|}
\hline
   {\bf R}$=$ 500 & & & & &&   \\ \hline
 {\bf T}$=$50,  $\gamma=0.6667$ & 1.0000 &    1.0000  &   1.0000  &   1.0000  &   0.9800  &   0.9200\\
   {\bf T}$=$50,  $\gamma=0.6500$  &1.0000   & 1.0000 &   0.9600 &   0.9960  &  0.8200  &  0.7800\\
    {\bf T}$=$50,  $\gamma=0.6400$ &  1.0000 &   1.0000  &  0.9600  &  0.9960  &  0.8200  &  0.7800\\
     {\bf T}$=$50,  $\gamma=0.6300$ &0.9920  &  1.0000  &  0.9920  &  0.9960 &   0.9600 &   0.0040\\
     \hline
      {\bf T}$=$100, $\gamma=0.6667$ &1.0000   &  1.0000  &   0.9880  &  1.0000   &  0.9200 &   0.9880\\
      {\bf T}$=$100, $\gamma=0.5900$ &  1.0000   &  1.0000  &   1.0000 &    1.0000 &    0.9960  &   0.7040\\
         {\bf T}$=$100, $\gamma=0.5500$ & 1.0000   & 1.0000   & 0.9880   & 1.0000   & 0.9560  &  0.9920\\
                      {\bf T}$=$100, $\gamma=0.5000$ & 1.0000  &  0.9960  &  0.9960  &  1.0000 &   0.9560  &  0.1360\\
    \hline
   {\bf T}$=$500, $\gamma=0.3704$ &1.0000 &   1.0000   & 1.0000  &  1.0000  &  1.0000  &  0.9520  \\
        {\bf T}$=$500, $\gamma=0.3077$ &   1.0000   &  1.0000   &  1.0000  &   1.0000 &    1.0000  &   0.9080 \\
         {\bf T}$=$500, $\gamma=0.2990$ &   1.0000  &    1.0000 &     1.0000  &    1.0000  &    0.9960 &     0.9080\\
           {\bf T}$=$500, $\gamma=   0.2650 $ &   1.0000&    1.0000 &   0.9920 &   1.0000  &  0.9720  &  0.6400\\
 \hline
 {\bf T}$=$1000, $\gamma=0.3704$ &1.0000  &    1.0000    &  1.0000   &   1.0000     &1.0000 &    1.0000     \\
  {\bf T}$=$1000, $\gamma=0.3077$ &  1.0000  &  1.0000 &   1.0000   & 0.9880  &  0.9520 &   0.3880\\
     {\bf T}$=$1000, $\gamma= 0.2650$ & 1.0000  &  1.0000  &  1.0000  &  1.0000 &   1.0000 &   0.9800\\
      {\bf T}$=$1000, $\gamma=0.2500$ &   1.0000  &  1.0000  &  1.0000  &  1.0000 &   0.9920 &   0.9480\\
   \hline
 \end{tabular}\label{T4}
\end{center}
\end{table}

 \begin{table}[!h]
\caption{\textbf{\emph{Example 2. Empirical power $M(T)=T^{2\gamma},$ $\beta =1/4,$   $\mathbf{k}_{n,j,h,l},$ $n=h=1,2,3,$ $\alpha =0.05$} }}
\begin{center}
\begin{tabular}{|l|llllll|}
\hline
   {\bf R}$=$ 500 & & & & &&   \\ \hline
 {\bf T}$=$50,  $\gamma=0.6667$ &     1.0000  &  1.0000  &  1.0000  &  1.0000  &  0.9800  &  0.9200\\
   {\bf T}$=$50,  $\gamma=0.6500$  &  0.9960   & 1.0000  &  0.9640   & 1.0000   & 0.8240  &  0.7800\\
     {\bf T}$=$50,  $\gamma=0.6400$ &  1.0000  &  1.0000 &   0.9600  &  0.9960 &   0.8160  &  0.7920\\
     {\bf T}$=$50,  $\gamma=0.6300$ & 0.9920 &   1.0000  &  0.9920  &  0.9960 &   0.9600 &   0.0040\\
     \hline
      {\bf T}$=$100, $\gamma=0.6667$ &  1.0000 &   1.0000   & 0.9880  &  1.0000  &  0.9200   & 0.9840\\
       {\bf T}$=$100, $\gamma=0.5900$ &  1.0000   & 1.0000   & 1.0000  &  1.0000  &  0.9960  &  0.7160\\
          {\bf T}$=$100, $\gamma=0.5500$ & 1.0000    & 1.0000   &  0.9880  &   1.0000  &   0.9560   &  0.9920\\
            {\bf T}$=$100, $\gamma=0.5000$ &  1.0000  &   0.9960  &   0.9960   &  1.0000   &  0.9560    & 0.1360\\
    \hline
   {\bf T}$=$500, $\gamma=0.3704$ & 1.0000 &   1.0000 &   1.0000   & 1.0000  &  0.9960 &   0.9560 \\
      {\bf T}$=$500, $\gamma=0.3077$ & 1.0000   & 1.0000   & 1.0000  &  1.0000  &  1.0000 &   0.9080\\
        {\bf T}$=$500, $\gamma=0.2800$ &  1.0000  &  1.0000 &   0.9920   & 1.0000  &  0.9720 &   0.6440\\
       {\bf T}$=$500, $\gamma=0.2650$ &   1.0000  &  1.0000  &  0.9920 &   1.0000  &  0.9720  &  0.6440\\
     \hline
 {\bf T}$=$1000, $\gamma=0.3704$ & 1.0000  &  1.0000  &  1.0000  &  1.0000  &  1.0000  & 1.0000    \\
  {\bf T}$=$1000, $\gamma=0.3077$ &    1.0000 &   1.0000  &  1.0000  &  0.9880  &  0.9520  &  0.3640\\
    {\bf T}$=$1000, $\gamma= 0.2650$ & 1.0000  &  1.0000  &  1.0000 &   1.0000  &  1.0000  &  0.9800\\
  {\bf T}$=$1000, $\gamma=0.2500$ &   1.0000  &  1.0000  &  1.0000  &  1.0000 &   0.9920 &   0.9480\\
   \hline
 \end{tabular}\label{T5}
\end{center}
\end{table}

 \begin{table}[!h]
\caption{\textbf{\emph{Example 3. Empirical power $M(T)=T^{2\gamma},$ $\beta =1/4,$   $\mathbf{k}_{n,j,h,l},$ $n=h=1,2,3,$ $\alpha =0.05$} }}
\begin{center}
\begin{tabular}{|l|llllll|}
\hline
   {\bf R}$=$ 500 & & & & &&   \\ \hline
 {\bf T}$=$50,  $\gamma=0.6667$ &   1.0000  &  1.0000  &  0.9960   & 1.0000  &  0.9800  &  0.9240\\
 {\bf T}$=$50,  $\gamma=0.6500$  & 1.0000 &   1.0000 &   0.9640 &   0.9960   & 0.8320  &  0.7800  \\
 {\bf T}$=$50,  $\gamma=0.6400$ &  1.0000  &  1.0000  &  0.9560  &  0.9960  &  0.8160  &  0.7880\\
     {\bf T}$=$50,  $\gamma=0.6300$ & 0.9920 &   1.0000 &   0.9920 &   0.9960  &  0.9600   & 0.0080\\
     \hline
      {\bf T}$=$100, $\gamma=0.6667$ &1.0000  &  1.0000  &  0.9880  &  1.0000 &   0.9200 &   0.9840\\
         {\bf T}$=$100, $\gamma=0.5900$ &  1.0000  &  1.0000  &  1.0000   & 1.0000 &   1.0000    & 0.7120\\
               {\bf T}$=$100, $\gamma=0.5500$ &   1.0000   & 1.0000   & 0.9880 &   1.0000   & 0.9560  &  0.9920\\
                 {\bf T}$=$100, $\gamma=0.5000$ &   1.0000 &   0.9960  &  0.9960  &  1.0000  &  0.9600  &  0.1440\\
    \hline
   {\bf T}$=$500, $\gamma=0.3704$ &  1.0000 &   1.0000 &   1.0000  &  1.0000  &  0.9960  &  0.9560    \\
    {\bf T}$=$500, $\gamma=0.3500$ &  1.0000  &  0.9680  &  0.9400   & 0.9960  &  0.9920   & 0.9920\\
        {\bf T}$=$500, $\gamma=0.3077$ &     1.0000   & 1.0000  &  0.9960  &  0.9920 &   0.7280  &  0.2960\\
                    {\bf T}$=$500, $\gamma= 0.2650$ &   1.0000  &  1.0000  &  0.9920 &   1.0000   & 0.9720   & 0.6480\\
     \hline
     {\bf T}$=$1000, $\gamma=0.3704$  & 1.0000  &  1.0000  &  1.0000  &  1.0000  &  1.0000  & 1.0000    \\
      {\bf T}$=$1000, $\gamma=0.3077$ & 1.0000  &   1.0000    & 1.0000   &  0.9880   &  0.9520   &  0.3880 \\
       {\bf T}$=$1000, $\gamma= 0.2650$ &   1.0000 &    1.0000   &  1.0000  &   1.0000 &    0.9960  &   0.9440\\
 {\bf T}$=$1000, $\gamma=0.2500$&  1.0000 &   1.0000  &  1.0000  &  1.0000  &  0.9880   & 0.9080  \\
   \hline
 \end{tabular}\label{T6}
\end{center}
\end{table}
\subsection{Discussion}
 The numerical results obtained in Section \ref{simuconsistency} illustrate Theorem \ref{th2} on the  consistency of the test. Particularly, one can observe that the asymptotic analysis ($T\to \infty$),  under sparse  spatial   scenarios,  requires lower values of sieve parameter $k(T) $ than the ones considered in Remark \ref{condiii}. That is the reason why we consider in this simulation study $k(T)=\mathcal{O}(\log (T)),$ $T\to \infty.$ This $T$--varying truncation scheme  is optimal for low model  complexity, when the functional data display smoothness.  Note that  our consistency analysis  requires the sieve parameter $k(T)$  ensures   the dominant subspace to be  included  into  one  of the  eigenspaces of the Laplace Beltrami operator generated by the elements of the sieve basis with dimension $k(T)$.

 Regarding the empirical size, one can observe, in Table \ref{T3bb},  that increasing the spherical  sampling frequency, a better performance is obtained for all the functional sample sizes analyzed. Note that, as expected,  for the sparsest spherical data scenario
($\gamma =0.25$), good approximations to the theoretical value $\alpha =0.05$ are also obtained, when the functional sample size increases, under low values of the sieve parameter.

The  effect of the spherical sampling frequency on the empirical power values has been illustrated in Tables \ref{T4}--\ref{T6}  for small functional sample sizes, where the minimum threshold for $\gamma ,$ in the definition of  $M(T),$ is intricate,  and strongly depends  on the spherical smoothness, and the memory of our functional time series model. In that sense, Tables \ref{T4}--\ref{T6} display some  critical lowest values of $\gamma ,$ which can be decreased  when the functional sample size  $T$ increases. As expected, a better performance of the proposed testing procedure is obtained when  increasing the number of observations  in time and/or space. Hence,  one can observe an increasing of the  empirical power values  at the six random projections, approaching the value one. Note that all the numerical results are obtained under  weak restrictions on smoothness since $s=3$  in condition (ii), i.e., the values of our  functional data set a.s.  lye in the Sobolev space $H^{3}(\mathbb{S}_{2}).$  Thus,  from the embedding theorems of Sobolev spaces into H\"older spaces, the functional data values a.s. lay in the  H\"older space $\mathcal{H}^{\rho }(\mathbb{S}_{2})$ of order $\rho < s-d/2=2$.

\section{Concluding remarks}
\label{cr}
Condition (\ref{cf}) is a key condition in our  strong consistency analysis of the test. In practice, getting a  balance between dimensionality in space and time constitutes a challenging topic where smoothness and memory of the underlying model play a crucial role.
Dimensionality in space is reflected  by the number of observable  spatial   random locations $M(T)=T^{\gamma },$ characterized by  parameter $\gamma$  in Remark \ref{condiii}. In our case, from Theorem \ref{th1}, low values of $\gamma $ also require low values of $\widetilde{\alpha } ,$ in the definition of
the $T$--varying sieve parameter $k(T)=T^{\widetilde{\alpha } },$ since $\widetilde{\alpha } <\gamma /2,$ under the formulation of  condition (iii)  in Remark \ref{condiii}. Specifically, from Theorem \ref{th2}, when  $\gamma -\widetilde{\alpha } <\frac{2\widetilde{\alpha }s}{d},$  spatial sparsity can be handled by increasing the number of temporal nodes (i.e., the sampling frequency in the temporal  spectral domain), under the memory condition $1/4<l_{\alpha }<1/2.$   Alternatively, when $\frac{2\widetilde{\alpha } s}{d} <\gamma -\widetilde{\alpha } $ holds, condition (\ref{cf}) is satisfied under  spatial smoothness, i.e., for $s$ sufficiently large, under the memory condition $1/4<l_{\alpha }<1/2.$ In this last case, the increasing of the topological dimension $d$ of the compact manifold also leads to stronger conditions on the local  regularity of the functional data under  sparse spatial observations.

In a near future, we will focus on  air pollution data analysis.  LRD spectral  testing  in spherical functional times series has substantial interest for this data set, collected across the globe or specific atmospheric  zones, using satellite measurements or discrete monitoring networks. They are mapped onto a sphere as a functional time series of continuous surface variation, or  spherical functional data. In these applications, air pollution levels (e.g., PM$_{2.5},$ PM$_{10}$ or Ozone) are treated  not as a collection of single points, but as a sequence of random functions with spherical support, i.e., a sequence of spherical random fields correlated in time. These observed spherical maps, generated from satellite instruments, and ground sensors do not provide a perfectly  continuous  snapshot of the Earth. Instead, data are available at discrete spatial locations with measurement errors, yielding the observation model  (\ref{reg}) introduced in Section \ref{dodata}. Note that air pollution exhibits strong persistence due to overarching climate dynamics, continuous background industrial emissions, and planetary-scale atmospheric transport systems. Thus, current pollution anomalies remain highly correlated, having important effects with those far in the future, requiring to go  beyond SRD framework, considered, for example,  in \cite{Alvarezlieb2019} for the analysis of such a data set.

\section*{Acknowledgements}
\noindent This work has been supported in part by projects PID2022--142900NB-I00, PID2024-158017NB-I00 and PID2020-116587GB-I00, financed by
\linebreak MCIU/AEI/10.13039/501100011033 and by FEDER UE,  and by CEX2020-001105-M/MCIN/AEI/ 10.13039/501100011033. It has also been partially supported  by regional grants ED431C 2021/24 and ED431C 2025/03 (Grupos competitivos),  financed by Xunta de Galicia through European Regional Development Funds (ERDF).

\section*{Appendix}
\begin{appendix}
\ref{SBLRDF} provides    the preliminary elements in the  functional spectral  domain, required in the formulation of the LRD test.
The proofs of the main  results of this paper are given in  \ref{supproofres}.

\setcounter{section}{0}
\section{Spectral background for  the LRD  test}
\label{SBLRDF}

Consider  the family of nuclear covariance operators $\left\{ \mathcal{R}_{\tau},\ \tau\in \mathbb{Z}\right\},$ with $\mathcal{R}_{\tau} =E[X_{s}\otimes X_{s+\tau}]=E[ X_{s+ \tau}\otimes X_{s}],$ for every $s,\tau \in \mathbb{Z}.$ For each $\tau\in \mathbb{Z},$
 $$\mathcal{R}_{\tau}(\phi)(\psi)=\int_{\mathbb{M}_{d}\times \mathbb{M}_{d}}r_{\tau}(x,y)\phi (x)\psi(y) d\nu (x)d\nu(y).$$

 The  Fourier transform in time $\{ \mathcal{F}_{\omega },\ \omega \in [-\pi,\pi]\}$ of the covariance operator family $\left\{ \mathcal{R}_{\tau},\ \tau\in \mathbb{Z}\right\},$ is obtained as the sum of the operator  series:
      \begin{eqnarray}& &
  \mathcal{F}_{\omega}
\underset{\mathcal{S}(L^{2}(\mathbb{M}_{d},d\nu; \mathbb{C}))}{=}  \frac{1}{2\pi} \sum_{\tau \in \mathbb{Z}}\exp\left(-i\omega \tau\right)\mathcal{R}_{\tau}\nonumber\\
 & &\underset{\mathcal{S}(L^{2}(\mathbb{M}_{d},d\nu; \mathbb{C}))}{=}\sum_{n\in \mathbb{N}_{0}}f_{n}(\omega)\sum_{j=1}^{\Gamma (n,d)}S_{n,j}^{d}\otimes \overline{S_{n,j}^{d}},\quad \omega \in [-\pi,\pi],
  \label{sdo2}
\end{eqnarray}
   \noindent where $\underset{\mathcal{S}(L^{2}(\mathbb{M}_{d},d\nu; \mathbb{C}))}{=}$ denotes the equality   in the  Hilbert--Schmidt operator norm $\|\cdot\|_{\mathcal{S}(L^{2}(\mathbb{M}_{d},d\nu; \mathbb{C}))},$ with $\mathcal{S}(L^{2}(\mathbb{M}_{d},d\nu; \mathbb{C}))$ being the space of Hilbert--Schmidt operators  on $L^{2}(\mathbb{M}_{d},d\nu; \mathbb{C}).$

   Note that the last equality in (\ref{sdo2}) holds under the assumption of
  invariance   of the covariance and spectral density operators of $X$  with respect to  the group of isometries of $\mathbb{M}_{d}.$
    The system $\{S_{n,j}^{d}, \ j=1,\dots, \Gamma (n,d),\ n\in \mathbb{N}_{0}\}$ denotes the   eigenfunctions  of the Laplace--Beltrami operator $\Delta_{d}$ on $L^{2}\left(\mathbb{M}_{d},d\nu , \mathbb{C}\right)$   (see, e.g., \cite{Gine75}; \cite{Helgason59}).
For every $n\in \mathbb{N}_{0},$ $\Gamma (n,d)$ represents the dimension of the eigenspace $\mathcal{H}_{n}$ associated with the eigenvalue $\lambda_{n}(\Delta_{d})$ of the Laplace Beltrami operator $\Delta_{d}$  (see, e.g., Section 2.1 in \cite{MaMalyarenko}).  The operator integrals are understood as improper operator Stieltjes integrals which converge strongly (see, e.g., Section 8.2.1 in  \cite{Ramm05}).

    The  fDFT $\widetilde{X}_{\omega }^{(T)}$  of $X$ provides an infinite-dimensional analogous to the discrete Fourier transform of the data. For a
  functional sample  \linebreak $\left\{ X_{t},\ t=0,\dots,T-1\right\}$ of size $T\geq 2,$ the  fDFT  is computed as
 \begin{eqnarray}\widetilde{X}_{\omega }^{(T)}(x)&=&\frac{1}{\sqrt{2\pi T}}\sum_{t=0}^{T-1}X_{t}(x)\exp(-i\omega t),\quad x\in \mathbb{M}_{d},\ \omega \in [-\pi,\pi].\label{fDFT}
\end{eqnarray}
\noindent    The empirical version of the cumulant operator of order two of the fDFT $\widetilde{X}_{\omega }^{(T)}$ over the diagonal $\omega \in [-\pi,\pi]$ is the  periodogram operator, denoted here as $\mathcal{P}^{(T)}_{\omega }.$ The  kernel  $p_{\omega }^{(T)}(x,y)$ of the periodogram operator  $\mathcal{P}^{(T)}_{\omega }=\widetilde{X}_{\omega}^{(T)}\otimes  \widetilde{X}^{(T)}_{-\omega}$  is given by,
 \begin{eqnarray}
  p_{\omega }^{(T)}(x,y)&=&\frac{1}{2\pi T}\sum_{t=0}^{T-1}\sum_{s=0}^{T-1}X_{t}(x)X_{s}(y)\exp(-i\omega [t-s]),\ \forall  x,y\in \mathbb{M}_{d},\nonumber\\ \label{periodogro}
\end{eqnarray}
\noindent    for every $\omega \in [-\pi,\pi].$ Thus, the   spectral cumulant operator $\mathcal{F}_{\omega }^{(T)}$ of order $2$ of the fDFT  $\widetilde{X}_{\omega }^{(T)}$  over the frequency diagonal has kernel $f_{\omega }^{(T)}(x,y),$ which can be computed as  $$f_{\omega }^{(T)}(x,y)=\mbox{cum}\left(\widetilde{X}^{(T)}_{\omega}(x), \widetilde{X}^{(T)}_{-\omega}(y)\right)=E\left[p^{(T)}_{\omega }(x,y)\right], \ x,y\in \mathbb{M}_{d}, \ \omega \in [-\pi,\pi].$$
Operator $\mathcal{F}_{\omega }^{(T)}$ can be written in terms of the  F\'ejer kernel  \begin{equation}F_{T}(\omega )=\frac{1}{T}\sum_{t=0}^{T-1}
\sum_{s=0}^{T-1}\exp\left(-i(t-s)\omega \right)=\frac{1}{T}\left[\frac{\sin\left(T\omega/2\right)}{\sin(\omega/2)}\right]^{2}
\label{fejkernel}
\end{equation}\noindent  as follows:
\begin{eqnarray}\mathcal{F}_{\omega }^{(T)}
&= &\left[F_{T}*\mathcal{F}_{\bullet
}\right](\omega )\nonumber\\ &&=
\int_{-\pi}^{\pi} F_{T}(\omega - \xi)
\mathcal{F}_{\xi} d\xi,\quad T\geq 2, \quad  \omega \in [-\pi,\pi].\label{convfejerkernel}\end{eqnarray}

As in the finite-dimensional case, one can consider a tapering in time of the functional data, leading to the introduction in the spectral domain of  the weighted periodogram operator  $\widehat{\mathcal{F}}_{\omega }^{(T)}$  (smoothed version in the   functional frequency domain of the periodogram operator)  having kernel $\widehat{f}_{\omega }^{(T)}(x,y),$ given,
  for each   $\omega \in [-\pi,\pi],$ by
 \begin{eqnarray}
&& \hspace*{-1cm} \widehat{f}_{\omega }^{(T)}(x,y)=\left[\frac{2\pi}{T}\right]\sum_{s=0}^{T-1}  W^{(T)}\left(\omega - \frac{2\pi s}{T} \right)
 p_{ \frac{2\pi  s}{T} }^{(T)}(x,y) ,\ x,y\in \mathbb{M}_{d},\label{enp}
\end{eqnarray}
\noindent where the weight function $W^{(T)}$  is such that
$W^{(T)}(x) = \sum_{j\in \mathbb{Z}}\frac{1}{B_{T}} W\left(\frac{x + 2\pi j}{B_{T}}\right),$
with $B_{T}$ being the positive  bandwidth parameter. Note that $W$ is  assumed to be a real--valued, positive, even, and bounded in variation    function defined on $\mathbb{R},$  with $W(x) =0$, if  $ |x|\geq 1,$
 $\int_{\mathbb{R}} \left|W(x)\right|^{2}dx <\infty,$ and   $\int_{\mathbb{R}} W(x)dx =1.$

The LRD scenario tested has been introduced  in   \cite{RuizMedina2022} in the spectral domain, under the following functional semiparametric modelling
 \begin{equation}
\mathcal{F}_{\omega }=\mathcal{M}_{\omega }|\omega|^{-\mathcal{A}},\quad \omega \in [-\pi,\pi],\label{LRD}
\end{equation}
\noindent where  the invariant positive self--adjoint operators  $\mathcal{M}_{\omega }$ and $|\omega|^{-\mathcal{A}}$  are composed  to define  $\mathcal{F}_{\omega }.$
 Note that, here, $\mathcal{A}$ denotes the LRD operator on  \linebreak $L^{2}(\mathbb{M}_{d},d\nu; \mathbb{C}),$ defining the functional parameter in the spectral semiparametric framework adopted.
  Operator $|\omega|^{-\mathcal{A}}$ in (\ref{LRD})  is interpreted as in the framework of operator self-similar processes (see, e.g., \cite{Charac14} and   \cite{Rackauskasv2}). In this framework,
  $\mathcal{A}$ would play the role of  operator--valued Hurst coefficient in the setting of  fractional Brownian motion. The regular  spectral operator   $\mathcal{M}_{\omega }$ is such that $X$ is Markovian  when the null space of $\mathcal{A}$ coincides with $L^{2}(\mathbb{M}_{d},d\nu; \mathbb{C}).$ Assume that $\mathcal{A}$ satisfies
\begin{equation}
 \mathcal{A}(f)(g)=\int_{\mathbb{M}_{d}\times \mathbb{M}_{d}}f(x)g(y)\sum_{n\in \mathbb{N}_{0}}\alpha (n)\sum_{j=1}^{\Gamma (n,d)}S_{n,j}^{d}(x) \overline{S_{n,j}^{d}}(y)d\nu(x)d\nu(y),
 \label{lerdop}
 \end{equation}

 \noindent  for $f,g \in C^{\infty }(\mathbb{M}_{d}),$
 the space of infinitely differentiable functions with compact support contained  in $\mathbb{M}_{d},$  where sequence $\left\{ \alpha (n),\ n\in \mathbb{N}_{0}\right\}$  is such that
$l_{\alpha }\leq \alpha (n)\leq L_{\alpha },$  for every $n\in \mathbb{N}_{0},$  and  $l_{\alpha }, L_{\alpha }\in (0,1/2).$
   Hence, $\mathcal{A}$ and $\mathcal{A}^{-1}$ are in the space $\mathcal{L}(L^{2}(\mathbb{M}_{d}, d\nu,\mathbb{C}))$ of bounded linear operators on $L^{2}(\mathbb{M}_{d}, d\nu,\mathbb{C}),$ with the norm  $\left\|\cdot \right\|_{\mathcal{L}(L^{2}(\mathbb{M}_{d}, d\nu,\mathbb{C}))}.$  In particular,
     $\left\|\mathcal{A}\right\|_{\mathcal{L}(L^{2}(\mathbb{M}_{d}, d\nu,\mathbb{C}))}<1/2.$
Note that   operator  $|\omega|^{-\mathcal{A}}$   is also  interpreted in the weak sense as  (\ref{lerdop}).
 For each $\omega \in [-\pi,\pi],$ operator  $\mathcal{M}_{\omega }$ in (\ref{LRD}) is  a  trace operator with kernel $\mathcal{K}_{\mathcal{M}_{\omega }}(x,y)$  satisfying
 $\mathcal{K}_{\mathcal{M}_{\omega }}(x,y)=\sum_{n\in \mathbb{N}_{0}}M_{n}(\omega )\sum_{j=1}^{\Gamma (n,d)}S_{n,j}^{d}\otimes \overline{S_{n,j}^{d}}(x, y),$ $x,y\in \mathbb{M}_{d},$ in terms of the sequence of positive eigenvalues
  $\{M_{n}(\omega ),\ n\in \mathbb{N}_{0}\}.$     For each $n\in \mathbb{N}_{0},$ $M_{n}(\omega ),$ $\omega \in [-\pi,\pi],$   is a  continuous positive slowly varying  function at
$\omega =0$ in the Zygmund's sense (see Definition 6.6 in \cite{Beran17}, and  Assumption IV in   \cite{RuizMedina2022}).
As commented, $X$ displays SRD, when $\alpha (n)=0,$ for every $n\in \mathbb{N}_{0}.$
 Note that, under the above introduced LRD scenario $\int_{-\pi}^{\pi}\left\|\mathcal{F}_{\omega }\right\|_{\mathcal{S}(L^{2}(\mathbb{M}_{d},d\nu, \mathbb{C}))}^{2}d\omega <\infty.$

\section{Proof of the results}
\label{supproofres}
\subsection{Proof of Theorem \ref{th1}}
\label{proofth1}
For each $\xi \in \Omega \backslash \Lambda_{0},$  with $\mathcal{P}(\Lambda_{0})=0,$ denote  for $t=0,\dots,T-1,$
 \begin{eqnarray}
&&\hspace*{-1cm}\widehat{\mathbf{X}}_{k(T)}(t,\xi)=\left(\widehat{X}_{00}(t,\xi),\  \dots, \dots, \widehat{X}_{N(k(T))\Gamma (N(k(T)),d)}(t,\xi)\right)^{\prime}\nonumber\\ &&\hspace*{2.5cm}=(B_{k(T),M(T)}^{\prime }B_{k(T),M(T)})^{-}B_{k(T),M(T)}^{\prime }\mathbf{X}_{t}\nonumber\\
&&\hspace*{3cm}+(B_{k(T),M(T)}^{\prime }B_{k(T),M(T)})^{-}B_{k(T),M(T)}^{\prime }\boldsymbol{\varepsilon}_{t},\nonumber\\
\label{epv}
\end{eqnarray}
\noindent whose entries $\widehat{X}_{nj}(t,\xi),$ $j=1,\dots,\Gamma (n,d),$ $n=0,\dots,N(k(T)),$ respectively approximate, from observation model (\ref{reg}), the entries
$X_{nj}(t,\xi),$  $j=1,\dots,$\linebreak $\Gamma (n,d),$ $n=0,\dots,N(k(T))$ of  vector $\mathbf{X}_{k(T)}(t,\xi),$  (see equations (\ref{eqkshexpansion}) and   (\ref{LSFC})).
Here,  \begin{eqnarray}&&\mathbf{X}_{t}=\left(X_{t}(z_{1}),\dots, X_{t}(z_{M(T)})\right)^{\prime },\quad
\boldsymbol{\varepsilon}_{t}=\left(\varepsilon_{1,t},\dots, \varepsilon_{M(T),t}\right)^{\prime }.\nonumber\end{eqnarray}

 We compute the mean  quadratic error
\begin{eqnarray}
&&\int_{\Omega \times \mathbb{M}_{d}}[X_{t}(z,\xi)-\widehat{\mathbf{X}}_{t,k(T),M(T)}(z)]^{2}d\nu(z) \mathcal{P}(d\xi)\nonumber\\
&&=\int_{\Omega \times \mathbb{M}_{d}}\left[\sum_{n=0}^{N(k(T))}\sum_{j=1}^{\Gamma (n,d)}[X_{nj}(t,\xi)-\widehat{X}_{nj}(t,\xi)]S_{n,j}(z)\right.\nonumber\\
&&\hspace*{1cm}\left.+
\sum_{n=N(k(T))+1}^{\infty }\sum_{j=1}^{\Gamma (n,d)}X_{nj}(t,\xi)S_{n,j}(z) \right]^{2}d\nu(z) \mathcal{P}(d\xi),\  t=0,\dots,T-1.\nonumber\\
\label{qe}
\end{eqnarray}

From  Parseval identity in equation (\ref{qe}), we obtain (see, e.g., \cite{Gine75})
\begin{eqnarray}
&&\int_{\Omega \times \mathbb{M}_{d}}[X_{t}(z,\xi)-\widehat{\mathbf{X}}_{t,k(T),M(T)}(z)]^{2}d\nu(z) \mathcal{P}(d\xi)\nonumber\\
&&=\sum_{n=0}^{N(k(T))}\sum_{j=1}^{\Gamma (n,d)}\int_{ \Omega}\left[X_{nj}(t,\xi)-\widehat{X}_{nj}(t,\xi)\right]^{2}\mathcal{P}(d\xi)\nonumber\\ &&
\hspace*{3cm}+
\sum_{n=N(k(T))+1}^{\infty }\sum_{j=1}^{\Gamma (n,d)}\int_{\Omega }X_{nj}^{2}(t,\xi)\mathcal{P}(d\xi)\nonumber\\
&&=\int_{\Omega }\left\|\mathbf{X}_{k(T)}(t,\xi)-\widehat{\mathbf{X}}_{k(T)}(t,\xi)\right\|^{2}\mathcal{P}(d\xi)\nonumber\\
&&\hspace*{3cm}+\int_{\Omega }
\sum_{n=N(k(T))+1}^{\infty }\sum_{j=1}^{\Gamma (n,d)}X_{nj}^{2}(t,\xi)\mathcal{P}(d\xi),
 \label{aaa}
\end{eqnarray}
\noindent  for $ t=0,\dots,T-1.$

Under (ii), almost surely (a.s.) in $\xi \in\Omega,$ the following asymptotic order holds:
\begin{equation}\sum_{n=N(k(T))+1}^{\infty }\sum_{j=1}^{\Gamma (n,d)}X_{nj}^{2}(t,\xi)=\mathcal{O}\left([k(T)]^{-2s/d}\right),\quad  T\to \infty.\label{l2em}
\end{equation}

In equation (\ref{aaa}), under  (i)--(iii), and (\ref{asidentity5}), Theorem 1 in \cite{Newey97}  can be applied considering  exponent $ s/d$ of the sieve parameter $k(T),$   and  $\widehat{\mathbf{X}}_{k(T)}$ being identified with the   parameter vector  in the  linear model (see also \cite{ChenChristensen15}). Thus, from this theorem, we obtain,  for $\xi\in \Omega \backslash  \Lambda_{0},$
\begin{equation}\left\|\mathbf{X}_{k(T)}(t,\xi)-\widehat{\mathbf{X}}_{k(T)}(t,\xi)\right\|^{2} =\mathcal{O}\left(\frac{k(T)}{M(T)}+[k(T)]^{-2s/d}\right),
\quad T\to \infty.\label{l2emb}\end{equation}

 From equations (\ref{l2em}) and (\ref{l2emb}), applying Dominated Convergence Theorem,  we have
\begin{eqnarray}
&&\hspace*{-2cm}\int_{\Omega \times \mathbb{M}_{d}}[X_{t}(z,\xi)-\widehat{\mathbf{X}}_{t,k(T),M(T)}(z)]^{2}d\nu(z) \mathcal{P}(d\xi)\nonumber\\
&&=\mathcal{O}\left(\frac{k(T)}{M(T)}+[k(T)]^{-2s/d}\right),\quad T\to \infty,
\label{eqmsecr}
\end{eqnarray}
 \noindent  for $t=0,\dots,T-1,$ as we wanted to prove.
 \subsection{Proof of Proposition \ref{pr1GH0}}
\label{pr1}
  Consider
  \begin{eqnarray}&&\widehat{\mathcal{S}_{B_{T}}}-E[\widehat{\mathcal{S}_{B_{T}}}]\nonumber\\ &&=[\widehat{\mathcal{S}_{B_{T}}}-\mathcal{S}_{B_{T}}]+[\mathcal{S}_{B_{T}}-E[\mathcal{S}_{B_{T}}]]+[E[\mathcal{S}_{B_{T}}]-
  E[\widehat{\mathcal{S}_{B_{T}}}]]\nonumber\\
  &&= S_{1}(T)+S_{2}(T)+S_{3}(T).
  \label{eqerrt}
  \end{eqnarray}

   Under SRD,  from equations (\ref{a1}), (\ref{a2}) and (\ref{a3}), assuming conditions in  Theorem \ref{th1}, and Theorem 2.2 in  \cite{Ruiz-MedinaCrujeiras24}, $$S_{1}(T)\leq  h_{8}(T)=\mathcal{O}\left(T^{1-\min\left\{(\gamma -\widetilde{\alpha})/2, \widetilde{\alpha} s/d\right\}}\right)+\mathcal{O}\left(B_{T}^{-1/2}T^{-1/2}\right),\ T\to \infty, $$ \noindent in the norm of the space
 $L^{2}(\Omega \times \mathbb{M}_{d}^{2},\mathcal{P}\otimes d\nu\otimes d\nu).$ Thus,  under the condition
 $\min\left\{(\gamma -\widetilde{\alpha})/2, \widetilde{\alpha} s/d\right\}>1,$ $S_{1}(T)\to 0,$ as $T\to \infty ,$ in the norm of the space
 $L^{2}(\Omega \times \mathbb{M}_{d}^{2},\mathcal{P}\otimes d\nu\otimes d\nu).$ As $T\to \infty ,$
 $S_{2}(T)$  converges in probability distribution to $Y_{0}^{(\infty)}$  from
 Theorem 2.2   in  \cite{Ruiz-MedinaCrujeiras24}.  The convergence to zero of $S_{3}(T)$ in the norm of the space
 $L^{2}( \mathbb{M}_{d}^{2}, d\nu\otimes d\nu)$ follows from the convergence to zero of $S_{1}(T)$  in the norm of the space
 $L^{2}(\Omega \times \mathbb{M}_{d}^{2},\mathcal{P}\otimes d\nu\otimes d\nu),$ by applying Jensen's inequality. Slutsky's Lemma then leads to the desired result.

\subsection{Proof of Lemma \ref{lem1}}
\label{lem1bb}
From triangle inequality,
\begin{eqnarray}
&&\left\|
\int_{-\pi}^{\pi}E_{H_{1}}\left[\widehat{\mathcal{P}}_{ \omega }^{(T)} \right]-\mathcal{F}_{\omega }d\omega d\omega \right\|_{L^{2}(\mathbb{M}_{d}^{2},d\nu \otimes d\nu, \mathbb{C})}\nonumber\\ &&
\leq \left\|\int_{-\pi}^{\pi}E_{H_{1}}\left[\widehat{\mathcal{P}}_{ \omega }^{(T)} -\mathcal{P}_{\omega}^{(T)}\right]d\omega \right\|_{L^{2}(\mathbb{M}_{d}^{2},d\nu \otimes d\nu, \mathbb{C})}\nonumber\\
&&+\left\|\int_{-\pi}^{\pi}E_{H_{1}}\left[\mathcal{P}_{\omega}^{(T)}\right]-\mathcal{F}_{\omega }d\omega \right\|_{L^{2}(\mathbb{M}_{d}^{2},d\nu \otimes d\nu, \mathbb{C})}=S_{1}(T)+S_{2}(T).\nonumber\\
 \label{e1l1}
\end{eqnarray}
Applying  Lemma 3.1 in \cite{Ruiz-MedinaCrujeiras24}, \begin{equation}S_{2}(T)=\mathcal{O}(T^{-1}),\  T\to \infty.\label{eqpf}
\end{equation}

From equation (\ref{a2}), and Jensen's inequality, applying  Dominated Convergence Theorem, we obtain
\begin{eqnarray}
&&
[S_{1}(T)]^{2}\leq h_{9}(T)=\mathcal{O}\left(T^{1-\min\left\{(\gamma -\widetilde{\alpha}), 2\widetilde{\alpha} s/d\right\}}\right).\label{eqk2}\end{eqnarray}

The result follows from equations (\ref{e1l1})--(\ref{eqk2}).

\subsection{Proof of Corollary \ref{lem3}}
\label{apcor3}
From triangle inequality,
\begin{eqnarray}
&&\left\|\int_{-\pi}^{\pi}E_{H_{1}}\left[\widetilde{\widehat{\mathcal{F}}}^{(T)}_{\omega }\right] -\int_{\mathbb{R}}W(\xi)\mathcal{F}_{\omega -\xi B_{T}}d\xi d\omega \right\|_{L^{2}(\mathbb{M}_{d}^{2},d\nu \otimes d\nu, \mathbb{C})}\nonumber\\
&&\leq \left\|\int_{-\pi}^{\pi}E_{H_{1}}\left[\widetilde{\widehat{\mathcal{F}}}^{(T)}_{\omega }-\widehat{\mathcal{F}}^{(T)}_{\omega }\right]d\omega\right\|_{L^{2}(\mathbb{M}_{d}^{2},d\nu \otimes d\nu, \mathbb{C})}\nonumber\\
&&+\left\|\int_{-\pi}^{\pi}E_{H_{1}}\left[\widehat{\mathcal{F}}^{(T)}_{\omega }\right]-
\int_{\mathbb{R}}W(\xi)\mathcal{F}_{\omega -\xi B_{T}}d\xi d\omega \right\|_{L^{2}(\mathbb{M}_{d}^{2},d\nu \otimes d\nu, \mathbb{C})}=S_{1}(T)+S_{2}(T).\nonumber\\
\label{eqcor1h}
\end{eqnarray}

Applying  Corollary 3.2 in \cite{Ruiz-MedinaCrujeiras24}, \begin{equation}S_{2}=\mathcal{O}(B_{T}^{-1}T^{-1})+\mathcal{O}(T^{-1}),\quad T\to \infty.\label{eqcr}
\end{equation}

Consider now the term $S_{1}^{2}(T).$ From  Jensen's inequality, we obtain
 \begin{eqnarray}
&&S_{1}^{2}(T)\leq \int_{-\pi}^{\pi}E_{H_{1}}\left[\left\|\widetilde{\widehat{\mathcal{F}}}^{(T)}_{\omega }-\widehat{\mathcal{F}}^{(T)}_{\omega }\right\|^{2}_{L^{2}(\mathbb{M}_{d}^{2}, d\nu \otimes d\nu, \mathbb{C})}\right]d\omega .\nonumber\\
\label{eqfi}\end{eqnarray}
   From  equation (\ref{a3}),  applying Dominated Convergence Theorem in (\ref{eqfi}), we  have
   \begin{equation}
 S_{1}^{2}(T)\leq h_{10}(T)=\mathcal{O}\left(T^{1-2\min\left\{(\gamma -\widetilde{\alpha} )/2, \widetilde{\alpha} s/d\right\}}B_{T}^{-1}\right),\ T\to \infty.
 \label{eqs11}
 \end{equation}

Equation (\ref{eqapphh}) then follows from equations (\ref{eqcor1h})--(\ref{eqs11}).

\subsection{Proof of Proposition  \ref{cor1}}
\label{pproposition2}
In what follows, $\widetilde{\mathcal{F}}^{(T)}_{\omega }$ denotes the integral operator with kernel
  \begin{eqnarray}&&\widetilde{f}^{(T)}_{\omega }(x,y)=E\left[\widehat{p}_{\omega }^{(T)}(x,y)\right]\nonumber\\ &&=\sum_{n\in \mathbb{N}_{0}}\sum_{j=1}^{\Gamma (n,d)}
  \widetilde{f}^{(T)}_{n}(\omega )S_{n,j}^{d}\otimes \overline{S_{n,j}^{d}}(x,y),\quad  x,y\in \mathbb{M}_{d},\label{feemp}
  \end{eqnarray}
   \noindent where, as before,  $\widehat{p}_{\omega }^{(T)}(x,y)$   is the kernel of  the plug--in periodogram operator $\widehat{\mathcal{P}}_{ \omega }^{(T)}$ in (\ref{pluperop}), and

   \begin{equation}\widetilde{f}^{(T)}_{n}(\xi )=\int_{-\pi}^{\pi} F_{T}(\xi - \omega)\widetilde{f}_{n}(\omega )d\omega,\quad \xi \in [-\pi,\pi],\ n\in \mathbb{N}_{0},\label{feempb}\end{equation}
\noindent with $F_{T}$ being  the F\'ejer kernel introduced in (\ref{fejkernel}), and, for $\omega \in [-\pi,\pi],$  $\widetilde{f}_{n}(\omega ),$ $n\in \mathbb{N}_{0},$ being
 such that
\begin{eqnarray}& &
  \widetilde{\mathcal{F}}_{\omega}
\underset{\mathcal{S}(L^{2}(\mathbb{M}_{d},d\nu; \mathbb{C}))}{=}  \frac{1}{2\pi} \sum_{\tau \in \mathbb{Z}}\exp\left(-i\omega \tau\right)\widehat{\mathcal{R}}_{\tau}=\sum_{n\in \mathbb{N}_{0}}\sum_{j=1}^{\Gamma (n,d)}\widetilde{f}_{n}(\omega )S_{n,j}^{d}\otimes \overline{S_{n,j}^{d}}\nonumber\\
&&\hspace*{3.5cm}\widehat{\mathcal{R}}_{\tau} =E[\widehat{X}_{s}\otimes \widehat{X}_{s+\tau}],\quad s,\tau \in \mathbb{Z}.\nonumber
\end{eqnarray}

From (\ref{feemp})--(\ref{feempb}),
\begin{eqnarray}
&&\left\| \int_{[-\sqrt{B_{T}}/2, \sqrt{B_{T}}/2]}E_{H_{1}}[\widetilde{\widehat{\mathcal{F}}}_{\omega  }^{(T)}]\frac{d\omega}{\sqrt{B_{T}}}\right\|_{\mathcal{S}(L^{2}(\mathbb{M}_{d},d\nu, \mathbb{C}))}
\nonumber\\
&&
 \geq
\left\| \int_{[-\sqrt{B_{T}}/2, \sqrt{B_{T}}/2]}E_{H_{1}}[\widetilde{\widehat{\mathcal{F}}}_{\omega  }^{(T)}]\frac{d\omega}{\sqrt{B_{T}}}\right\|_{\mathcal{L}(L^{2}(\mathbb{M}_{d},d\nu, \mathbb{C}))}
\nonumber\\
&&
= \sup_{n\in \mathbb{N}_{0}}\left|\int_{-\pi}^{\pi}\frac{1}{B_{T}}W\left(\frac{\xi }{B_{T}}\right)
\widetilde{f}^{(T)}_{n}(\xi)d\xi \right|\nonumber\\
&&\hspace*{2cm}
\geq g(T)=\mathcal{O}(B_{T}^{-1/2-l_{\alpha }}),\quad  T\to \infty,
\label{eqaltdiv}
\end{eqnarray}\noindent   where $\mathcal{L}(L^{2}(\mathbb{M}_{d},d\nu, \mathbb{C}))$ is the space of bounded linear operators on $L^{2}(\mathbb{M}_{d},d\nu, \mathbb{C}),$ and $\widetilde{f}^{(T)}_{n}(\xi)$ satisfies (\ref{feempb}).

\subsection{Proof of Theorem \ref{consisimse}}
\label{consist}
From triangle, Cauchy--Schwartz,  and Jensen's inequalities,
\begin{eqnarray}
&&
\int_{-\pi}^{\pi}E_{H_{1}}\left\|\widetilde{\widehat{\mathcal{F}}}_{\omega }^{(T)}-E_{H_{1}}[\widetilde{\widehat{\mathcal{F}}}_{\omega }^{(T)}]\right\|_{\mathcal{S}(L^{2}(\mathbb{M}_{d},d\nu, \mathbb{C}))}^{2}d\omega \nonumber\\
&&\leq 6\int_{-\pi}^{\pi}E_{H_{1}}\left\|\widetilde{\widehat{\mathcal{F}}}_{\omega }^{(T)}-\widehat{\mathcal{F}}_{\omega }^{(T)}
\right\|_{L^{2}(\mathbb{M}_{d}^{2},d\nu\otimes d\nu, \mathbb{C})}^{2}d\omega \nonumber\\
&&+2\int_{-\pi}^{\pi}E_{H_{1}}\left\|\widehat{\mathcal{F}}_{\omega }^{(T)}-E_{H_{1}}\left[\widehat{\mathcal{F}}_{\omega }^{(T)}\right]
\right\|_{L^{2}(\mathbb{M}_{d}^{2},d\nu\otimes d\nu, \mathbb{C})}^{2}d\omega \nonumber\\ &&=S_{1}(T)+S_{2}(T).
\label{thconstproof}
\end{eqnarray}
Considering  now equation (\ref{a3}),   Dominated Convergence Theorem leads to
\begin{eqnarray}
&& S_{1}(T)\leq  h_{11}(T)=\mathcal{O}(B_{T}^{-1}T^{-1})+\mathcal{O}\left(T^{(1+\beta)-2\min\left\{(\gamma -\widetilde{\alpha } )/2  , \widetilde{\alpha } s/d\right\}}\right),\ T\to \infty.
\nonumber\\
\label{dctheqth2}
\end{eqnarray}
Finally, from  Theorem 4.2 in  \cite{Ruiz-MedinaCrujeiras24},
\begin{eqnarray}
&& S_{2}(T)\leq h_{12}(T)=\mathcal{O}(B_{T}^{-1}T^{-1}),\quad \ T\to \infty.\label{dctheqth2b}
\end{eqnarray}

From equations (\ref{thconstproof})-- (\ref{dctheqth2b}) we obtain (\ref{ca2}).

\subsection{Proof of Theorem \ref{th2}}
\label{maresult}
The proof follows from Proposition   \ref{cor1}  and Theorem \ref{consisimse}, adopting  a similar methodology to the proof of Theorem 4.5 in \cite{Ruiz-MedinaCrujeiras24}. Specifically,
the plug--in test statistic operator $\widehat{\mathcal{S}}_{B_{T}}$ is  reformulated as
\begin{eqnarray}
&&\hspace*{-1cm}\widehat{\mathcal{S}}_{B_{T}}=\sqrt{B_{T}T}\int_{[-\sqrt{B_{T}}/2, \sqrt{B_{T}}/2]}E_{H_{1}}\left[\widetilde{\widehat{\mathcal{F}}}_{\omega }^{(T)} \right]\frac{d\omega}{\sqrt{B_{T}}}
\nonumber\\
&& \circ
\left[\mathbb{I}_{L^{2}(\mathbb{M}_{d},d\nu, \mathbb{C})}+\left[\int_{[-\sqrt{B_{T}}/2, \sqrt{B_{T}}/2]}\left(\widetilde{\widehat{\mathcal{F}}}_{\omega }^{(T)}-
E_{H_{1}}\left[\widetilde{\widehat{\mathcal{F}}}_{\omega}^{(T)}\right]\right)\frac{d\omega}{\sqrt{B_{T}}}  \right]\right.
\nonumber\\
&&
\left.\circ \left[\int_{[-\sqrt{B_{T}}/2, \sqrt{B_{T}}/2]}E_{H_{1}}\left[\widetilde{\widehat{\mathcal{F}}}_{\omega }^{(T)} \right]\frac{d\omega}{\sqrt{B_{T}}} \right]^{-1}\right],
\label{st}\end{eqnarray}
\noindent where $\circ$ means the composition of operators,
$\mathbb{I}_{L^{2}(\mathbb{M}_{d},d\nu, \mathbb{C})}$ denotes  the identity operator on the space
 $L^{2}(\mathbb{M}_{d},d\nu, \mathbb{C}),$  and $\left[\int_{[-\sqrt{B_{T}}/2, \sqrt{B_{T}}/2]}E_{H_{1}}\left[\widetilde{\widehat{\mathcal{F}}}_{\omega}^{(T)}\right]\frac{d\omega}{\sqrt{B_{T}}} \right]^{-1}$ is  the inverse of operator
  $\int_{[-\sqrt{B_{T}}/2, \sqrt{B_{T}}/2]}E_{H_{1}}\left[\widetilde{\widehat{\mathcal{F}}}_{\omega}^{(T)}\right]\frac{d\omega}{\sqrt{B_{T}}}.$

  From Proposition    \ref{cor1},  as $T\to \infty,$
\begin{eqnarray}&&\left\|\sqrt{B_{T}T}\int_{[-\sqrt{B_{T}}/2, \sqrt{B_{T}}/2]}E_{H_{1}}\left[\widetilde{\widehat{\mathcal{F}}}_{\omega }^{(T)}\right]\frac{d\omega}{\sqrt{B_{T}}} \right\|_{\mathcal{S}(L^{2}(\mathbb{M}_{d},d\nu, \mathbb{C}))}
\nonumber\\
&& \hspace*{4cm}\geq g(T)=\mathcal{O}\left(T^{1/2}B_{T}^{-l_{\alpha }}\right).\nonumber\\
\label{eqcor1}
\end{eqnarray}

 The following inequality holds:
\begin{eqnarray}
&&E_{H_{1}}\left\|\int_{[-\sqrt{B_{T}}/2, \sqrt{B_{T}}/2]}\left[\widetilde{\widehat{\mathcal{F}}}_{\omega }^{(T)}-
E_{H_{1}}\left[\widetilde{\widehat{\mathcal{F}}}_{\omega }^{(T)}\right]\right]\frac{d\omega}{\sqrt{B_{T}}}
\right.
\nonumber\\
&&
\hspace*{2cm}\circ \left.\left[\int_{[-\sqrt{B_{T}}/2, \sqrt{B_{T}}/2]}E_{H_{1}}\left[\widetilde{\widehat{\mathcal{F}}}_{\omega}^{(T)}\right]\frac{d\omega}{\sqrt{B_{T}}}  \right]^{-1}
\right\|^{2}_{\mathcal{S}(L^{2}(\mathbb{M}_{d},d\nu, \mathbb{C}))}\nonumber\\
&&\leq \left\|
\left[\int_{[-\sqrt{B_{T}}/2, \sqrt{B_{T}}/2]}E_{H_{1}}\left[\widetilde{\widehat{\mathcal{F}}}_{\omega}^{(T)}\right]\frac{d\omega}{\sqrt{B_{T}}} \right]^{-1} \right\|^{2}_{\mathcal{L}(L^{2}(\mathbb{M}_{d},d\nu, \mathbb{C}))}\nonumber\\
&&\hspace*{0.5cm} \times
E_{H_{1}}\left\|\int_{[-\sqrt{B_{T}}/2, \sqrt{B_{T}}/2]}\left[\widetilde{\widehat{\mathcal{F}}}_{\omega}^{(T)}-
E_{H_{1}}\left[\widetilde{\widehat{\mathcal{F}}}_{\omega}^{(T)}\right]\right]\frac{d\omega}{\sqrt{B_{T}}} \right\|^{2}_{\mathcal{S}(L^{2}(\mathbb{M}_{d},d\nu, \mathbb{C}))}.
\label{eqseconsumnum}
\end{eqnarray}
Again, applying Proposition    \ref{cor1},  as $T\to \infty,$\begin{eqnarray}&&\hspace*{-0.3cm} \left\|\left[\int_{[-\sqrt{B_{T}}/2, \sqrt{B_{T}}/2]}E_{H_{1}}\left[\widetilde{\widehat{\mathcal{F}}}_{\omega}^{(T)}\right]\frac{d\omega}{\sqrt{B_{T}}}\right]^{-1}
\right\|^{2}_{\mathcal{L}(L^{2}(\mathbb{M}_{d},d\nu, \mathbb{C}))}\leq h_{13}(T)=\mathcal{O}\left(B_{T}^{2l_{\alpha }+1}\right).\label{rest1}\nonumber\\
\end{eqnarray}
 \noindent   From  Jensen inequality,
\begin{eqnarray}&&
E_{H_{1}}\left\|\int_{[-\sqrt{B_{T}}/2, \sqrt{B_{T}}/2]}\left[\widetilde{\widehat{\mathcal{F}}}_{\omega}^{(T)}-
E_{H_{1}}\left[\widetilde{\widehat{\mathcal{F}}}_{\omega}^{(T)}\right]\right]\frac{d\omega}{\sqrt{B_{T}}} \right\|^{2}_{\mathcal{S}(L^{2}(\mathbb{M}_{d},d\nu, \mathbb{C}))}\nonumber\\
&&\hspace*{-0.7cm} \leq  \int_{[-\sqrt{B_{T}}/2, \sqrt{B_{T}}/2]}E_{H_{1}}\left\|\widetilde{\widehat{\mathcal{F}}}_{\omega}^{(T)}-
E_{H_{1}}\left[\widetilde{\widehat{\mathcal{F}}}_{\omega}^{(T)}\right]\right\|^{2}_{\mathcal{S}(L^{2}(\mathbb{M}_{d},d\nu, \mathbb{C}))}\frac{d\omega}{\sqrt{B_{T}}}.\label{rest2bb}
\end{eqnarray}
 \noindent Hence,  from equations (\ref{eqseconsumnum})--(\ref{rest2bb}), applying   Theorem \ref{consisimse},   as $T\to \infty,$ \begin{eqnarray}&&\hspace*{-0.5cm} E_{H_{1}}\left\|\int_{[-\sqrt{B_{T}}/2, \sqrt{B_{T}}/2]}\left[\widetilde{\widehat{\mathcal{F}}}_{\omega}^{(T)}-
E_{H_{1}}\left[\widetilde{\widehat{\mathcal{F}}}_{\omega}^{(T)}\right]\right]\frac{d\omega}{\sqrt{B_{T}}} \right.\nonumber\\
&&\left.\hspace*{0.5cm}\circ \left[\int_{[-\sqrt{B_{T}}/2, \sqrt{B_{T}}/2]}E_{H_{1}}\left[\widetilde{\widehat{\mathcal{F}}}_{\omega}^{(T)}\right]\frac{d\omega}{\sqrt{B_{T}}}  \right]^{-1}
\right\|^{2}_{\mathcal{S}(L^{2}(\mathbb{M}_{d},d\nu, \mathbb{C}))}
\nonumber\\
&&\leq h_{14}(T)=\mathcal{O}\left(T^{-1}B_{T}^{2l_{\alpha }-1/2}\right)+\mathcal{O}\left(T^{-\beta (2l_{\alpha }+1)+1+3\beta /2-\min\left\{\gamma -\widetilde{\alpha }  , 2\widetilde{\alpha } s/d\right\}}\right).\nonumber\\
\label{sicv}
\end{eqnarray}
 From  equation (\ref{sicv}), applying Chebyshev's inequality,
 \begin{eqnarray}
 &&
 P\left[\left\|\int_{[-\sqrt{B_{T}}/2, \sqrt{B_{T}}/2]}\left[\widetilde{\widehat{\mathcal{F}}}_{\omega}^{(T)}-
E_{H_{1}}\left[\widetilde{\widehat{\mathcal{F}}}_{\omega}^{(T)}\right]\right]\frac{d\omega}{\sqrt{B_{T}}} \right.\right.\nonumber\\
&&\left.\left.\hspace*{0.5cm}\circ \left[\int_{[-\sqrt{B_{T}}/2, \sqrt{B_{T}}/2]}E_{H_{1}}\left[\widetilde{\widehat{\mathcal{F}}}_{\omega}^{(T)}\right]\frac{d\omega}{\sqrt{B_{T}}}  \right]^{-1}\right\|_{\mathcal{S}(L^{2}(\mathbb{M}_{d},d\nu, \mathbb{C}))}>\varepsilon\right]\nonumber\\
&&\leq E_{H_{1}}\left\|\int_{[-\sqrt{B_{T}}/2, \sqrt{B_{T}}/2]}\left[\widetilde{\widehat{\mathcal{F}}}_{\omega}^{(T)}-
E_{H_{1}}\left[\widetilde{\widehat{\mathcal{F}}}_{\omega}^{(T)}\right]\right]\frac{d\omega}{\sqrt{B_{T}}} \right.\nonumber\\
&&\left.\hspace*{0.5cm}\circ \left[\int_{[-\sqrt{B_{T}}/2, \sqrt{B_{T}}/2]}E_{H_{1}}\left[\widetilde{\widehat{\mathcal{F}}}_{\omega}^{(T)}\right]\frac{d\omega}{\sqrt{B_{T}}}  \right]^{-1}
\right\|^{2}_{\mathcal{S}(L^{2}(\mathbb{M}_{d},d\nu, \mathbb{C}))}/\varepsilon^{2}\nonumber\\
&&\hspace*{-0.3cm}\leq h_{15}(T)/\varepsilon^{2}=\mathcal{O}(T^{-1}B_{T}^{2l_{\alpha }-1/2})+\mathcal{O}\left(T^{-\beta (2l_{\alpha }+1)+1+3\beta /2 -\min\left\{\gamma -\widetilde{\alpha }  , 2\widetilde{\alpha }s/d\right\}}\right).\nonumber\\
\label{feq}
\end{eqnarray}
Since   $l_{\alpha }>1/4,$  hence, $2l_{\alpha }-1/2=\rho>0,$   and,  for $B_{T}=T^{-\beta },$
$T^{-1}B_{T}^{2l_{\alpha }-1/2}=T^{-1-\beta \rho},$ with $\beta \in (0,1),$ and $\rho \in (0,1/2).$ Thus,  from condition (\ref{cf}) and  equation  (\ref{feq}),  Borel--Cantelli  lemma   leads,  as $T\to \infty,$ to
\begin{eqnarray}
&&\left\|\int_{[-\sqrt{B_{T}}/2, \sqrt{B_{T}}/2]}\left[\widetilde{\widehat{\mathcal{F}}}_{\omega}^{(T)}-
E_{H_{1}}\left[\widetilde{\widehat{\mathcal{F}}}_{\omega}^{(T)}\right]\right]\frac{d\omega}{\sqrt{B_{T}}}\right. \nonumber\\
&&\hspace*{2cm}\left. \circ \left[\int_{[-\sqrt{B_{T}}/2, \sqrt{B_{T}}/2]}E_{H_{1}}\left[\widetilde{\widehat{\mathcal{F}}}_{\omega}^{(T)}\right]\frac{d\omega}{\sqrt{B_{T}}}  \right]^{-1}\right\|_{\mathcal{S}(L^{2}(\mathbb{M}_{d},d\nu, \mathbb{C}))}\to_{a.s.} 0.\nonumber\\
\label{feqfv}
\end{eqnarray}
\noindent
  The a.s. divergence of $\left\|\widehat{\mathcal{S}}_{B_{T}}\right\|_{\mathcal{S}(L^{2}(\mathbb{M}_{d},d\nu, \mathbb{C}))},$ as  $T\to \infty ,$ follows   from   equations  (\ref{st}), (\ref{eqcor1}) and (\ref{feqfv}).

\end{appendix}
\end{document}